\documentclass[a4paper,english, 11pt]{article}
\usepackage[a4paper, margin=1.3in
]{geometry}
\usepackage[utf8]{inputenc}
\usepackage[T1]{fontenc}
\usepackage{lmodern}
\usepackage[english]{babel}
\usepackage{hyperref}

\usepackage{amsmath}
\usepackage{amsthm}
\usepackage{amsfonts}
\usepackage{amssymb}
\usepackage{mathtools}
\usepackage{siunitx}
\usepackage[linesnumbered, ruled, noend]{algorithm2e}
\usepackage{graphicx}
\usepackage{subcaption}
\usepackage{longtable}
\usepackage{multirow}
\usepackage{enumerate}
\usepackage{xcolor}
\usepackage{xspace}
\usepackage{tikz}
\usepackage[font=small,labelfont=bf]{caption}
\usepackage{array}
\usepackage{setspace}

\usetikzlibrary{fit}
\usetikzlibrary{shapes}
\usetikzlibrary{shapes.geometric}
\usetikzlibrary{backgrounds}
\usetikzlibrary{arrows}
\usetikzlibrary{arrows.meta}
\tikzset{
	-Latex,auto,node distance =0.5 cm and 0.5 cm,semithick,
	square/.style={regular polygon,regular polygon sides=4}
}

\usetikzlibrary{decorations}
\usetikzlibrary{decorations.text}
\usetikzlibrary{topaths}
\usetikzlibrary{patterns}
\usetikzlibrary{calc}
\usetikzlibrary{intersections}
\usetikzlibrary{through}
\usetikzlibrary{spy}
\usetikzlibrary{matrix}
\usetikzlibrary{lindenmayersystems}
\usetikzlibrary{external}
\usetikzlibrary{positioning}
\usetikzlibrary{3d}

\theoremstyle{plain}
\newtheorem{theorem}{Theorem}[section]

\newtheorem{proposition}[theorem]{Proposition}

\newtheorem{observation}[theorem]{Observation}

\usepackage[pagewise]{lineno}

\usepackage[
backend=biber,
style=apa,
urldate=short,
]{biblatex}

\AtNextBibliography{\small}

\DefineBibliographyStrings{english}{%
  urlseen = {Accessed on}, 
}

\DeclareFieldFormat{urldate}{(\bibstring{urlseen}\space#1)}
\DeclareFieldFormat{url}{\url{#1}}

\renewbibmacro*{url+urldate}{%
  \usebibmacro{url}
  \usebibmacro{urldate}%
}

\newcommand{\python}{{\sc python}\xspace}
\newcommand{\gurobi}{{\sc gurobi}\xspace}
\newcommand{\gurobiVersion}[2]{{\sc gurobi}~\oldstylenums{#1.#2}\xspace}

\newcommand{\cplexVersion}[2]{{\sc cplex}~\oldstylenums{#1.#2}\xspace}

\newcommand{\scipjack}{{\sc \mbox{scip-Jack}}\xspace}

\newcommand{\scipVersion}[3]{{\sc scip}~\oldstylenums{#1.#2.#3}\xspace}

\newcommand{\Z}{\mathbb{Z}}
\newcommand{\R}{\mathbb{R}}

\newcommand{\generalProblem}{\ensuremath{P}\xspace}

\newcommand{\optSolution}[1]{\ensuremath{v(#1)}\xspace}
\newcommand{\lpRelaxation}[1]{\ensuremath{LP_{#1}}\xspace}
\newcommand{\arcVariable}{\ensuremath{x}\xspace}
\newcommand{\nodeVariable}{\ensuremath{y}\xspace}

\newcommand{\quotaInstance}{\ensuremath{I_{\mathrm{Q}}}\xspace}
\newcommand{\quotaInstanceTrans}{\ensuremath{I_{\mathrm{QT}}}\xspace}

\newcommand{\terminals}{\ensuremath{T}\xspace}
\newcommand{\fixTerminals}{\ensuremath{\terminals_f}\xspace}
\newcommand{\fixTerminalsTrans}{\ensuremath{\terminals^{\prime}_f}\xspace}
\newcommand{\potTerminals}{\ensuremath{\terminals_p}\xspace}

\newcommand{\quota}{\ensuremath{Q}\xspace}
\newcommand{\cost}{\ensuremath{c}\xspace}
\newcommand{\scenicness}{\ensuremath{s}\xspace}
\newcommand{\scenicnessVertex}{\ensuremath{s^\mathrm{\vertex}}\xspace}

\newcommand{\setVertices}{\ensuremath{V}\xspace}

\newcommand{\vertex}{\ensuremath{v}\xspace}
\newcommand{\vertexNumbered}[1]{\ensuremath{\vertex_{#1}}\xspace}
\newcommand{\vertexCosts}{\ensuremath{w}\xspace}
\newcommand{\vertexProfits}{\ensuremath{q}\xspace}
\newcommand{\vertexCostsNumbered}[1]{\ensuremath{\vertexCosts_{#1}}\xspace}
\newcommand{\vertexProfitsNumbered}[1]{\ensuremath{\vertexProfits_{#1}}\xspace}

\newcommand{\setEdges}{\ensuremath{E}\xspace}
\newcommand{\edge}{\ensuremath{e}\xspace}

\newcommand{\edgeCost}[1]{\ensuremath{\cost_{#1}}\xspace}
\newcommand{\edgeCostVertices}[2]{\ensuremath{\cost_{#1#2}}\xspace}

\newcommand{\fullQSTPInstance}{\ensuremath{\quotaInstance = (\setVertices, \setArcs, \fixTerminals, \potTerminals, c, \vertexProfits, \quota)}\xspace}
\newcommand{\fullTransQSTPInstance}{\ensuremath{\quotaInstanceTrans = (\setVertices^{\prime}, \setArcs^{\prime}, \terminals^{\prime}, c^{\prime}, \vertexProfits^{\prime}, \quota)}\xspace}

\newcommand{\graph}{\ensuremath{G}\xspace}
\newcommand{\graphFull}{\ensuremath{\graph = (\setVertices,\setEdges)}\xspace}

\newcommand{\setArcs}{\ensuremath{A}\xspace}
\newcommand{\arc}{\ensuremath{a}\xspace}

\newcommand{\arcCostVertices}[2]{\ensuremath{\cost(#1,#2)}\xspace}
\newcommand{\diGraphFull}{\ensuremath{\ensuremath{D} = (\setVertices,\setArcs)}\xspace}
\newcommand{\incomingArcs}[1]{\ensuremath{\delta^-(#1)}\xspace}
\newcommand{\outgoingArcs}[1]{\ensuremath{\delta^+(#1)}\xspace}

\newcommand{\arcFlow}{\ensuremath{f}\xspace}

\usepackage{csquotes}
\usepackage{zibtitlepage}
\usepackage{authblk}

\makeatletter
\newcommand\blfootnote[1]{%
  \begingroup
  \renewcommand{\@makefntext}[1]{\noindent\makebox[1.8em][r]#1}
  \renewcommand\thefootnote{}\footnote{#1}%
  \addtocounter{footnote}{-1}%
  \endgroup
}
\makeatother

\begin{document}
\ZTPAuthor{
    Jaap Pedersen, Jann Michael Weinand, Chloi Syranidou, Daniel Rehfeldt
}

\ZTPInfo{The content of this report is also available as publication in European Journal of Operational Research. Please always cite as:

\fullcite{pedersen2024a}
}
\ZTPNumber{23-10}
\ZTPMonth{May}
\ZTPYear{2024}

\title{An efficient solver for large-scale onshore wind farm siting including cable routing}
\date{}

\author[a]{Jaap Pedersen\footnote{\ZTPOrcid{0000-0003-4047-0042}, corresponding author, pedersen@zib.de}$^,$}
\author[b]{Jann Michael Weinand\footnote{\ZTPOrcid{0000-0003-2948-876X}}$^,$}
\author[b]{Chloi Syranidou\footnote{\ZTPOrcid{0000-0002-3332-6635}}$^,$}
\author[a]{Daniel~Rehfeldt\footnote{\ZTPOrcid{0000-0002-2877-074X}}$^,$}

\affil[a]{Zuse Institute Berlin, Berlin, Germany}
\affil[b]{Forschungszentrum Jülich GmbH, Institute of Energy and Climate Research – Techno-economic Systems Analysis (IEK-3), 52425 Jülich, Germany}

\hypersetup{pdftitle={
                    An efficient solver for large-scale onshore wind farm siting including cable routing},
            pdfauthor={Jaap Pedersen, Jann Michael Weinand, Chloi Syranidou, Daniel Rehfeldt}}



\maketitle
\vspace{-20pt}
\begin{abstract}
    Existing planning approaches for onshore wind farm siting and grid integration often do not meet minimum cost solutions or social and environmental considerations. In this paper, we develop an exact approach for the integrated layout and cable routing problem of onshore wind farm planning using the Quota Steiner tree problem. Applying a novel transformation on a known directed cut formulation, reduction techniques, and heuristics, we design an exact solver that makes large problem instances solvable and outperforms generic MIP solvers. In selected regions of Germany, the trade-offs between minimizing costs and landscape impact of onshore wind farm siting are investigated. Although our case studies show large trade-offs between the objective criteria of cost and landscape impact, small burdens on one criterion can significantly improve the other criteria. In addition, we demonstrate that contrary to many approaches for exclusive turbine siting, grid integration must be simultaneously optimized to avoid excessive costs or landscape impacts in the course of a wind farm project. Our novel problem formulation and the developed solver can assist planners in decision-making and help optimize wind farms in large regions in the future.\blfootnote{\newline 
    
    The content of this report is also available as publication in European Journal of Operational Research. Please always cite as:

    \fullcite{pedersen2024a}}
\end{abstract}
{\small
\subsection*{Keywords}
OR in energy, combinatorial optimization, quota Steiner tree problem, graph theory 
}

\section{Introduction}
The deployment of low-carbon technologies is vital in order to mitigate climate change. As part of the transformation of the global energy system, low-cost wind energy has become an established source of electricity \parencite{Veers.2019}. Between 2000 and 2019, global wind turbine capacity increased by more than 20\% annually \parencite{Pryor.2020}, reaching 730 GW in 2020 \parencite{OurWorldinData.2021}, with a further increase of 50\% anticipated by the end of 2023 \parencite{Pryor.2020}. Experts also predict a sharp decline in the already low cost of wind energy by 2050 \parencite{Wiser.2016,Wiser.2021,Jansen.2020}. 

In recent years, however, the expansion of onshore wind has stalled in some countries and regions \parencite{IEA.2021}. Despite general approval, local interest groups, such as neighboring residents or communities, increasingly oppose the construction of onshore wind turbines \parencite{Reusswig.2016,Rand.2017,Weinand.2021}, especially if they are not involved in the planning process \parencite{Fast.2016,Boudet.2019}. One of the main reasons for this opposition is the visual impact of these installations on the landscape \parencite{Petrova.2016,Suskevics.2019,Wolsink.2018,Molnarova.2012,Spielhofer.2021}. This opposition is most pronounced against placement of wind turbines in landscapes of high aesthetic/scenic quality, whereas wind turbines placed in less attractive landscapes are more likely to be accepted \parencite{Molnarova.2012}. The impact of resistance is particularly evident in Germany, the country with the third-largest onshore wind capacity \parencite{Statista.2020} (around 58 GW in 2022 \parencite{CleanEnergyWire.2023}) and the fourth largest share of onshore wind in electricity generation worldwide (around 26\%, \parencite{Statista.2020b}). After record years in 2014 and 2017, with capacity additions of 4.8 GW and 5.3 GW, respectively, only 1.0 GW, 1.4 GW, 1.6 GW and 2.4 GW of new capacity was added in 2019, 2020, 2021 and 2022 \parencite{CleanEnergyWire.2021, CleanEnergyWire.2023}, respectively. The rapid expansion and development of onshore wind turbines has triggered an increase in local protest movements and lawsuits across the country \parencite{CleanEnergyWire.2019,Buck.2019}. Together with the hurdles erected by legislators for new wind turbines, this raises doubts as to whether the government's expansion target of an additional minimum of 57 GW by 2030 is feasible \parencite{CleanEnergyWire.2023}. 

In the case of larger onshore wind farms, overhead lines are usually used for grid connection, which has also led to social opposition due to the accompanying modification of the landscape \parencite{Reusswig.2016,Bertsch.2016}. Even the typically used underground cables can have a negative impact on the landscape, e.g. via the cutting of paths and protective strips in forests \parencite{Roth.2021c}. Furthermore, historical wind energy projects have made it apparent that grid operators often implement suboptimal grid connection plans. Although grid operators are obliged under the German Renewable Energy Sources Act to establish the most economically-favorable connection points, this is rarely the case in reality \parencite{BWE.2021}. A planning approach for turbine location and grid connection planning that would optimally take into account the central target criteria of cost efficiency and landscape impact could accelerate onshore wind expansion again and so support the achievement of expansion targets. 

\subsection{Literature on wind farm planning}

Recent research studies have dealt with the optimal siting of onshore wind turbines using multi-criteria objectives: in \textcite{Weinand.29.06.2021}, \textcite{Lehmann.2021b}, \textcite{Tafarte.2021}, optimal onshore wind sites were determined for the entirety of Germany on the basis of the turbine levelized cost of electricity in €/MWh, aestetic qualities of landscapes or disamenities for the inhabitants living close to them. In \textcite{Weinand.historic}, these analyses were extended to turbine siting across Europe, with a simulation rather than an optimization approach being taken due to the large number of turbines and scenarios. In these large-scale studies, grid integration is typically neglected due to the high complexity of the resulting combinatorial problem \parencite{Weinand.combinatorial}. According to the results of the grid integration heuristic for onshore wind turbines discussed in 
\textcite{McKenna.2021b}, however, the total costs are doubled on average when grid integration is taken into account. It is therefore imperative to include grid integration. Due to the complexity of the problem, grid connection planning is mostly formulated as a minimum spanning tree (MST) problem (e.g. \textcite{Wedzik.2016}), is mostly solved heuristically \parencite{Wu.2020,HerbertAcero.2014,Cazzaro.2020,Cazzaro.2022}, and turbine locations are specified as fixed (e.g. \textcite{Hertz.2017,Žarković.2017}).

In the design process of a wind farm the two main challenges are a) the wind farm layout optimization (WFLO) and b) the wind farm cable routing (WFCR). The WFLO maximizes the energy output by reducing wake effects. Each installed turbine introduces a \textit{wake}, a region of slower wind speed, resulting in lower energy yield of turbines downwind. The WFCR minimizes the cost to connect the chosen turbines to the power grid. Due to the complexity of both problems, they are usually solved sequentially. However, there are some promising examples for the design and cable routing of offshore wind farms that integrate these problems \parencite{Fischetti.2018,Fischetti.2021a, cazzaro2023}, but typically a fixed number of turbines are determined prior to the cable routing. Some of these studies also include some artificial (Steiner) nodes to provide for more flexibility for the routing problem. This is sufficient for a single wind farm in an offshore environment, in which there are fewer obstacles than in the case of onshore wind farm planning. In \textcite{fischetti2022}, an integrated layout and cable routing problem is discussed in the context of offshore wind farms. The problem is modeled as a mixed-integer linear program, which is then improved by using cutting techniques. Even though the model can handle a variable set of wind turbines, it was decided to fix the number, as done in practice of single wind farm design. All latter articles \parencite{Fischetti.2018, Fischetti.2021a, fischetti2022} discuss the design of a single offshore wind farm, solving a model including relevant technical details on graphs with 50~-~100 nodes. Recently, \textcite{cazzaro2023} also discuss the integrated layout and cable optimization problem. The authors approach the problem heuristically using an adapted Variable Neighborhood Search. Although the number of available positions is high, the cable routing is performed on a much smaller set of fixed turbine positions. The mentioned studies show the complexity of the WFLO, WFCR, and their integrated variant even for a single offshore wind farm. The problem becomes even more complex if the installation process and the maintenance of the wind farm are considered, which are both more cost intensive than in the case of an onshore wind farm, see, e.g., \textcite{amorosi2024, adeirawan2023, gutierrez-alcoba2019}. Due to different terrain characteristics, varying wind profiles over large distances, and legal regulations positioning of onshore wind turbines is much more involving than offshore. The goal of this study is to provide an approach for onshore wind planning on a large-scale. Therefore, the focus lies on the combinatorial aspect of the problem, i.e., which turbines to choose and where to route the cables.

\subsection{Literature, Algorithms and Software for Steiner Trees}

An efficient way to solve the problem for larger instances is to reduce the level of technical detail, e.g., incorporating the wake effect a priori, ignoring cable losses, or capacity constraints on cables and substations, etc. In this study, the problem of choosing a subset of possible wind turbines, including cable routing, is modeled as a variant of the Steiner tree problem in graphs (STP), a generalization of the MST problem. The STP is a classic NP-hard problem \parencite{karp1972}, and one of the most studied problems in combinatorial optimization \parencite{ljubic2021}. Given an undirected graph with non-negative edge costs, the STP aims to find a tree that interconnects a given set of special points (referred to as terminals) at minimum total cost. The STP and its variations arise in many real-world applications like network design problems in telecommunication, electricity, or in district heating, as well as other fields such as biology; see, e.g., \textcite{leitner2014}, \textcite{bolukbasi2018}, \textcite{ljubic2006}, and \textcite{klimm2020}, respectively. A vast amount of literature concerning STP and related problems exists and for a comprehensive overview of the topic, the reader is referred to a recent survey \parencite{ljubic2021}. In the context of wind farm design, the STP appears for example in the PhD thesis by \textcite{ridremont2019}, in which a robust cable network is designed; again, the method is only applied on networks with up to 100 nodes and 300 edges.

In general, network design problems consist of two parts: first, a subset of profitable customers must be selected, and second, a network must be designed to connect all chosen customers with the least cost. The trade-off between maximizing profits and minimizing costs can be modeled as a generalized version of the STP: the Prize-Collecting Steiner Tree Problem in graphs (PCSTP). However, in the transition towards a carbon-neutral energy system only focusing on this trade-off is insufficient. For example, network operators must assure security of demand, there are expansion targets to be met, or a region strives to cover its demand by locally operated renewable energy sources. These additional constraints can be modeled in a generalized version of the PCSTP, namely the quota PCSTP (QSTP), with the objective of minimizing costs while a minimum amount of "profits" is collected. 

Although a wide range of literature discusses the PCSTP from both the theoretical as well as the practical point of view (see the surveys in \textcite{costa2006} and \textcite{ljubic2021}), few studies have addressed the QSTP. The QSTP was initially formulated by \textcite{johnson2000}, who observed that it is a generalized version of the k minimum spanning tree (k-MST) problem, in which, given an edge-weighted graph G, a minimum costs subtree of G containing at least k vertices is constructed. The authors also propose a heuristical approach to solve the QSTP by introducing an increasing profit-multiplier $\alpha$, and solving a series of PCSTP instances. In each iteration, a new instance is constructed by multiplying the original profits by the increased $\alpha$ until the quota is fulfilled. This approach results in a trade-off curve that shows which quota can be collected at which prize. Drawbacks of this approach include the fact that multiple problems must be solved and it might not be clear how to choose $\alpha$, especially if costs and profits are not comparable, e.g., costs vs. energy potential; finally, the desired quota might not be captured in the trade-off curve. \textcite{haouari2006} propose a hybrid Lagrangian genetic algorithm to compute the lower and upper bounds for QSTP instances with up to 5000 edges. A robust version of the QSTP was introduced in \textcite{alvarez-miranda2013}. However, the authors only discuss a branch-and-cut approach for the robust PCSTP and its budget-constrained variant and not for the QSTP. Although a number of studies on exact solution approaches for the PCSTP have been conducted so far (e.g., \cite{ljubic2006, leitner2018, rehfeldt2022}), to the best of our knowledge no exact solution approach for the QSTP has been suggested in the literature. 

As mentioned above, we want to solve the integrated layout and cable routing problem on a large scale involving a high number of nodes, i.e., potential wind turbine positions and possibly Steiner nodes. The specialized Steiner tree solver \scipjack has proven that it can efficiently solve STP-related problems of magnitudes larger than what general out-of-the-box MIP solvers can even load into memory. \scipjack is backed by a number of theoretical and practical results, see, e.g., \textcite{rehfeldt2019, rehfeldt2022, RehfeldtKoch2023}. \scipjack is an extension of the non-commercial general MIP solver SCIP \parencite{BestuzhevaEtal2021OO}. Alongside various problem-specific reduction techniques and heuristics, \scipjack uses a branch-and-cut procedure to handle the exponential number of constraints induced by the Steiner cut-like constraints. These constraints will be introduced in Section \ref{sec:DCutForm}. For a detailed description of \scipjack's operational principle the reader is referred to \textcite{rehfeldt2021a}. By introducing the QSTP to \scipjack, we aim to solve the integrated layout and cable routing problem, although simplified, but on a much larger scale than the studies presented in the previous section.

\subsection{Contribution and structure}

In this study, we develop a planning instrument based on the Steiner tree approach, which is applicable to any region in Germany. This research is interdisciplinary and combines mathematical, landscape planning and energy management methods. In contrast to studies on offshore wind farm planning \parencite{Fischetti.2018, Fischetti.2021a, fischetti2022}, our approach aims to solve onshore wind turbine placements and cable routing on a regional level, resulting in much larger instances as a larger number of turbines are considered and more Steiner nodes are needed to assure flexibility in the cable routing. We apply our approach for the two cases of minimizing costs and the scenic impact of the installed turbines and cable connections. In Section \ref{sec:QSTP}, we first present the QSTP as a directed cut formulation. Then, we introduce a new transformation for the problem and prove the equivalence with the original formulation. We integrate the transformed QSTP into the exact Steiner tree solver \scipjack \parencite{RehfeldtKoch2023}. We further implement a shortest-path-based reduction technique and a primal heuristic for the QSTP. Subsequently, we apply our new methodological approaches in Section \ref{sec:Results} to demonstrate the significant performance improvements compared to standard solvers and to demonstrate the trade-offs between cost and landscape impacts in onshore wind farm planning for several German regions. Thereby, we highlight what previous (national) planning approaches must improve, before discussing and concluding our methods and results in Section \ref{sec:Discussion}.

\section{Quota Steiner tree problem}\label{sec:QSTP}

We model the integrated wind farm layout and cable routing problem as the Quota Steiner tree problem (QSTP). The available wind turbine positions are defined by the discrete set $\potTerminals$, which we call the set of potential terminals. Next, we call the discrete set of the grid's substations in the problem the set of fixed terminals $\fixTerminals$. There might be additional nodes in the graph called Steiner nodes, whose only purpose is to allow a more flexible cable routing. Finally, we have given a discrete set of edges E which represents the set of cable connection we can choose from. 

The QSTP is defined similarly to \textcite{johnson2000} with the addition of vertex costs: 
let \graphFull be an undirected graph. It is given a set of fixed terminals $\fixTerminals \subset \setVertices$ and a set of potential terminals $\potTerminals \subset \setVertices$ with $\fixTerminals \cap \potTerminals = \emptyset$. A node $\vertex\in V\setminus(\fixTerminals\cup\potTerminals)$ is called Steiner node. Each potential terminal $\vertex\in \potTerminals$ is associated with costs $\vertexCosts : \potTerminals\rightarrow \R_{>0}$ and quota profits $\vertexProfits: \potTerminals \rightarrow \R_{>0}$. In terms of wind farm planning, these quota profits represent the turbine's annual energy yield. The costs of each edge $(i,j) \in \setEdges$ is defined by $\cost : \setEdges \rightarrow \R_{\ge0}$. The goal is to find a tree $S = (\setVertices^\prime, \setEdges ^\prime) \subseteq \graph$ that contains all fixed terminals $\fixTerminals$ such that the total cost, i.e., cable and turbine costs:
\begin{equation}
    C(S) = \sum_{(i,j)\in \setEdges^\prime} \edgeCostVertices{ij} + \sum_{i \in \potTerminals \cap \setVertices^\prime} \vertexCostsNumbered{i}
\end{equation}
is minimized and a given quota $\quota\in \R_{>0}$ is fulfilled, i.e.:
\begin{equation}
    \quota(S) = \sum_{i \in \potTerminals \cap \setVertices^\prime} \vertexProfitsNumbered{i} \ge \quota \label{eq:quota_definition}
\end{equation}

In this study, the focus is on wind turbine siting and cable routing on a regional scale. Therefore, we assume that the quota represents, for example, a wind turbine expansion target. However, by setting the quota profit $\vertexProfitsNumbered{\vertex} = 1$ for all $\vertex\in\potTerminals$ and $\quota=k$ with $k\in \Z$, we can achieve that a minimum number of wind turbines have to be installed to fulfill \eqref{eq:quota_definition}.

\subsection{Directed cut formulation}\label{sec:DCutForm}
The QSTP is modeled as an Steiner arborescence problem (SAP) in an integer linear program (IP). For the general SAP see e.g. \textcite{wong1984}. The original undirected graph $\graph$ is transformed into a directed graph $\diGraphFull$ where $\setArcs\coloneqq \lbrace (i,j),(j,i) | \forall (i,j) \in \setEdges\rbrace$. By applying the idea of shifting the costs of a vertex $v$ onto the costs of its incoming arcs\footnote{Since the optimal solution is an arborescence, each vertex has only one incoming arc, so the vertex' costs are account for exactly ones if the vertex is included.} (see \textcite{ljubic2006}), 
the arc costs $\cost : \setArcs \rightarrow \R_{\ge0}$ are defined as:
\begin{align}
& \arcCostVertices{i}{j} = \begin{dcases} 
        \edgeCost{e} + \vertexCostsNumbered{j} & \text{if } j \in \potTerminals,\\
        \edgeCost{e}& \text{otherwise}
        \end{dcases} & \forall \arc=(i,j) \in \setArcs
\end{align}
where $\edgeCost{e}$ represents the cost of the corresponding edge $\edge=(i,j)$ in the original undirected graph $\graph$. For a subset of nodes $W\subset V$, we denote $\delta^+(W) = \lbrace (i,j) \in \setArcs: i \in W, j \in \setVertices\backslash W \rbrace$ as the set of outgoing arcs and $\delta^-(W) = \lbrace (i,j) \in A: i \in \setVertices\backslash W, j \in W \rbrace$ as the set of incoming arcs. For a single vertex $\vertexNumbered{i}$, we write $\delta^+(\lbrace \vertexNumbered{i}\rbrace) = \delta^+(\vertexNumbered{i})$ and $\delta^-(\lbrace \vertexNumbered{i}\rbrace) = \delta^-(\vertexNumbered{i})$. For any set $K$, we define $x(K) = \sum_{i \in K} x_i$. 

We introduce a binary variable $\arcVariable_{ij}$ for each $(i,j) \in \setArcs$ if the arc $(i,j)$ is contained in the Steiner tree ($\arcVariable_{ij} = 1$) or not ($\arcVariable_{ij} = 0$). Furthermore, let $\nodeVariable_k$ be a binary variable for each $k \in \potTerminals$ indicating whether or not the potential terminal $k$ is chosen. The rooted directed cut integer programming formulation (IP) of the QSTP with an arbitrary root $r \in \fixTerminals$ is given as follows:
\begin{align}
    \min \quad&c^T x &\label{eq:obj_qstp}\\
    \text{s.t.} \quad&&\nonumber\\
    & x(\delta^-(W)) \ge 1& \forall W \subset \setVertices, r\notin W, |W \cap \fixTerminals| \ge 1\label{Steiner cut}\\
    & x(\delta^-(W)) \ge y_i& \forall W \subset \setVertices, r\notin W, |W \cap \potTerminals| \ge 1, \vertexNumbered{i} \in \potTerminals \label{Steiner cut pot}\\
    & \sum_{\vertexNumbered{i} \in \potTerminals} \vertexProfitsNumbered{i} y_i \ge \quota & \label{quota cons}\\
    & x_{ij}, y_{k}  \in \lbrace 0,1 \rbrace &\forall (i,j) \in \setArcs, \forall \vertexNumbered{k} \in \potTerminals \label{eq:var_qstp}
\end{align}
The Steiner cut constraint \eqref{Steiner cut} requires that any subset of nodes of the graph $W \subset \setVertices$ that does not contain the root ($r\notin W$), but at least one fixed terminal, i.e., $|W \cap \potTerminals| \ge 1$, has at least one incoming arc. This guarantees that there exists a path from the root $r$ to each $t\in \fixTerminals$. The same goes for \eqref{Steiner cut pot}, i.e., a path exists from $r$ to the potential terminal $i\in \potTerminals$ if it is chosen to contribute to the quota constraint \eqref{quota cons} ($y_i = 1$). Let \fullQSTPInstance denote an instance of the QSTP. Note a simple observation for the QSTP:

\begin{observation}\label{obs:qstp_stp}
Considering an instance \fullQSTPInstance of the QSTP, if
\begin{align}
    \quota > \sum_{\vertexNumbered{i} \in \potTerminals} \vertexProfitsNumbered{i} - \min_{\vertexNumbered{i} \in \potTerminals}\vertexProfitsNumbered{i} \label{eq:quota_max}
\end{align}
then the QSTP reduces to the STP by defining the set of terminals $\terminals = \fixTerminals \cup \potTerminals$.
\end{observation}

In the case of \eqref{eq:quota_max}, one would need to choose all potential terminals, i.e., $y_i = 1$ for all $\vertexNumbered{i} \in \potTerminals$, to fulfill the quota constraint \eqref{quota cons}. However, by defining the set of fixed terminals by $\overline{\fixTerminals}= \fixTerminals\cup\potTerminals$ and $\overline{\potTerminals} = \emptyset$, as all potential terminal have to be chosen anyway, constraints \eqref{Steiner cut pot} and \eqref{quota cons} become obsolete. Thus, only \eqref{Steiner cut} remain and the problem is reduced to the original directed cut formulation for the STP which can be solved efficiently by specialized solvers like \scipjack. So, the QSTP generalizes the STP and, thus, the QSTP is also NP-hard. 

\subsubsection{Transformation}
One of the essential features of the existing, state-of-the-art \scipjack framework is the cut separation algorithm for the SAP. This algorithm makes it possible to solve LPs including the, exponentially many, constraints~\eqref{Steiner cut} in polynomial time. Initially, one starts with only a small subset of these constraints, and solves the LP(-relaxation). Next, the separation algorithm of \scipjack can be used to either find violated constraints~\eqref{Steiner cut} and add them to the LP, or to verify that all constraints are satisfied. If violated constraints have been added, the procedure is repeated. Otherwise, an optimal LP solution has been found. By the famous equivalence of optimization and separation, see \cite{grotschel1981ellipsoid}, this procedure is guaranteed to solve the LP in polynomial time (since violated constraints can be found in polynomial time). Because of the sophisticated separation algorithm of \scipjack, this procedure works very well in practice.
However, the, exponentially many, constraints \eqref{Steiner cut pot} cannot be separated by this algorithm; at least not without modifications. Therefore, we present a transformation of the problem formulation in the following that gets rid of the constraints \eqref{Steiner cut pot} and therefore allows us to directly apply the efficient separation algorithm of \scipjack.

For each potential terminal $t_i \in \potTerminals$, a new fixed terminal $t_{i^\prime}$ with a profit of $\vertexProfitsNumbered{t_{i^\prime}} = \vertexProfitsNumbered{t_i}$ is added. Furthermore, for each newly-added terminal $t_{i^\prime}$, an arc $(r,t_{i^\prime})$ with costs $\arcCostVertices{r}{t_{i^\prime}} = 0$ and an arc $(t_i, t_{i^\prime})$ with costs $\arcCostVertices{t_i}{t_{i^\prime}} = 0$ are added. Finally, each original potential terminal $t_i \in \potTerminals$ now becomes a Steiner node $i$ in the transformed graph. Let $\fixTerminalsTrans$ be the set of the newly-added terminals and $\terminals^{\prime} = \fixTerminals\cup \fixTerminalsTrans$ the set of all terminals in the transformed graph. 

Now, if we tried to keep the structure of the quota constraint \eqref{quota cons}, the quota profit $\vertexProfitsNumbered{t_{i^\prime}}$ would be collected if the newly-added fixed terminal $t_{i^\prime}$ is reached by the arc $(t_i, t_{i^\prime})$. However, by adding a constraint of form $b^T x \ge d$ with $b\in \R_{\ge0}^{A}$ and $b \ne \mathbf{0}$ and $d\in\R$ 
to the general SAP formulation, we can no longer guarantee the connectivity of the solution, which will be explained in the following. Let us consider only the Steiner cut constraints \eqref{Steiner cut}. As \textcite{goemans1993} 
point out, the convex hull of all $\arcVariable \in \Z_{\ge0}^{A}$ satisfying these constraints is of blocking type, i.e., its recession cone consists of all non-negative vectors. This means any solution feasible for \eqref{Steiner cut} can be extended by increasing $\arcVariable_\arc$ for any $\arc \in \setArcs$ and still be feasible. Only if we minimize over non-negative arc costs, an optimal solution yields a Steiner tree. However, by introducing $b^T x \ge d$
, one can add arcs (increasing the values of $\arcVariable$) to the solution until the ``$\ge$'' condition holds without violating the Steiner cut constraint \eqref{Steiner cut}. These arcs, however, do not necessarily have to be connected to the root component. To avoid this problem, instead of considering which potential terminals are taken into account, we are interested in which potential terminals are not considered, i.e., which newly-added terminal is directly connected to the root $r$ by the newly-added arc $(r,t_{i^\prime})$ and not via the original potential terminal. Thus, to fulfill the quota, there is an upper limit on how many newly-added terminals are connected directly to the root $r$; see \eqref{eq:quota_cons_trans}. 

\begin{figure}
    \begin{subfigure}{0.4\linewidth}
        \centering
        \begin{tikzpicture}[scale=1.5, every node/.style={scale=1.}]
            \tikzstyle{terminal} = [draw, fill, thick, minimum size=0.4, inner sep=2.6pt, color=red];
            \tikzstyle{pot_terminal} = [draw, fill, thick, minimum size=0.4, inner sep=2.6pt, color=blue];
            \tikzstyle{pot_terminald} = [draw, fill, densely dotted, minimum size=0.4, inner sep=2.6pt, color=blue!20];
            \tikzstyle{steiner} = [circle, fill, draw, thick, minimum size=0.2, inner sep=1.5pt];
            \tikzstyle{steinerd} = [circle, densely dotted, draw, minimum size=0.2, inner sep=1.5pt];
            \def\y{16}
            \node[terminal] (yt1) at (0.5+\y, .5) {};
            \node[terminal] (yt2) at (0.5+\y, 3.5) {} node[above of=yt2,yshift=-20pt]{$r$};
            \node[pot_terminal] (ypt1) at(-.5 + \y, 1.) {} node[left of=ypt1, xshift=10pt]{\scriptsize q=10};
            \node[pot_terminal] (ypt2) at(1.5 + \y, 1.2) {} node[right of=ypt2, xshift=-10pt]{\scriptsize q=40};
            \node[pot_terminal] (ypt3) at(0. + \y, 2.5) {} node[left of=ypt3, xshift=10pt]{\scriptsize q=20};
            \node[steiner] (ys1) at (0.5 + \y, 1.3) {};
            \node[steiner] (ys2) at (1.3 + \y, 2.3) {};
                
            \draw[thick, -latex](yt1) edge[bend left=20] node [midway, left=-2.0pt] {} (ys1); 
            \draw[thick, latex-](yt1) edge[bend left=-20] node [midway, right=-2.0pt] {} (ys1); 
            
            \draw[thick, -latex](ys1) edge[bend left=20] node [midway, left=-2.0pt] {} (ys2); 
            \draw[thick, latex-](ys1) edge[bend left=-20] node [midway, right=-2.0pt] {} (ys2); 
            \draw[thick, -latex](ys1) edge[bend left=20] node [midway, below=-2.0pt] {} (ypt1); 
            \draw[thick, latex-](ys1) edge[bend left=-20] node [midway, above=-2.0pt] {} (ypt1); 
            \draw[thick, -latex](ys1) edge[bend left=20] node [midway, left=-2.0pt] {} (ypt3); 
            \draw[thick, latex-](ys1) edge[bend left=-20] node [midway, right=-2.0pt] {} (ypt3); 
            
            \draw[thick, -latex](ys1) edge[bend left=20] node [midway, above=-2.0pt] {} (ypt2); 
            \draw[thick, latex-](ys1) edge[bend left=-20] node [midway, below=-2.0pt] {} (ypt2); 
            \draw[thick, -latex](ys2) edge[bend left=20] node [midway, above=-2.0pt] {} (ypt3); 
            \draw[thick, latex-](ys2) edge[bend left=-20] node [midway, above=-2.0pt] {} (ypt3); 
            \draw[thick, -latex](yt2) edge[bend left=20] node [midway, right=-2.0pt] {} (ypt3); 
            \draw[thick, latex-](yt2) edge[bend left=-20] node [midway, left=-2.0pt] {} (ypt3); 
            
            \draw[thick, -latex](yt2) edge[bend left=20] node [midway, above=-2.0pt] {} (ys2); 
            \draw[thick, latex-](yt2) edge[bend left=-20] node [midway, above=-2.0pt] {} (ys2); 
        \end{tikzpicture}
        \caption{Original instance.}
    \end{subfigure}
    \begin{subfigure}{0.6\linewidth}
    \centering
        \begin{tikzpicture}[scale=1.5, every node/.style={scale=1.}]
            \tikzstyle{terminal} = [draw, fill, thick, minimum size=0.4, inner sep=2.6pt, color=red];
            \tikzstyle{pot_terminal} = [draw, fill, thick, minimum size=0.4, inner sep=2.6pt, color=blue];
            \tikzstyle{pot_terminald} = [draw, fill, densely dotted, minimum size=0.4, inner sep=2.6pt, color=blue!20];
            \tikzstyle{steiner} = [circle, fill, draw, thick, minimum size=0.2, inner sep=1.5pt];
            \tikzstyle{steinerd} = [circle, densely dotted, draw, minimum size=0.2, inner sep=1.5pt];
            \def\y{16}
            \node[terminal] (yt1) at (0.5+\y, .5) {};
            \node[terminal] (yt2) at (0.5+\y, 3.5) {} node[above of=yt2,yshift=-20pt]{$r$};
            \node[steiner] (ypt1) at(-.5 + \y, 1.) {};
            \node[steiner] (ypt2) at(1.5 + \y, 1.2) {};
            \node[steiner] (ypt3) at(0. + \y, 2.5) {};
            \node[steiner] (ys1) at (0.5 + \y, 1.3) {};
            \node[steiner] (ys2) at (1.3 + \y, 2.3) {};
            \draw[thick, -latex](yt1) edge[bend left=20]   (ys1); 
            \draw[thick, latex-](yt1) edge[bend left=-20]  (ys1); 
            \draw[thick, -latex](ys1) edge[bend left=20]   (ys2); 
            \draw[thick, latex-](ys1) edge[bend left=-20]  (ys2); 
            \draw[thick, -latex](ys1) edge[bend left=20]   (ypt1);
            \draw[thick, latex-](ys1) edge[bend left=-20]  (ypt1);
            \draw[thick, -latex](ys1) edge[bend left=20]   (ypt3);
            \draw[thick, latex-](ys1) edge[bend left=-20]  (ypt3);
            \draw[thick, -latex](ys1) edge[bend left=20]   (ypt2);
            \draw[thick, latex-](ys1) edge[bend left=-20]  (ypt2);
            \draw[thick, -latex](ys2) edge[bend left=20]   (ypt3);
            \draw[thick, latex-](ys2) edge[bend left=-20]  (ypt3);
            \draw[thick, -latex](yt2) edge[bend left=20]   (ypt3);
            \draw[thick, latex-](yt2) edge[bend left=-20]  (ypt3);
            \draw[thick, -latex](yt2) edge[bend left=20]   (ys2); 
            \draw[thick, latex-](yt2) edge[bend left=-20]  (ys2); 
            \draw[thick, latex-](yt2) edge[bend left=-20]  (ys2); 
            \node[terminal](yt11) [left of=ypt1, xshift=-5pt] {} node[below left of=yt11, xshift=10pt, yshift=10pt]{\scriptsize $q^\prime=10$};
            \node[terminal](yt22) [right of=ypt2, xshift=5pt] {} node[below right of=yt22, xshift=-5pt, yshift=10pt]{\scriptsize $q^\prime=40$};
            \node[terminal](yt33) [left of=ypt3, xshift=-5pt] {} node[below left of=yt33, xshift=12pt, yshift=10pt]{\scriptsize $q^\prime=20$};
            \draw[thick, -latex](yt2) edge[bend left=-60] node [midway, left=2pt]{\scriptsize \textbf{0~}}  (yt11);
            \draw[thick, -latex](ypt1) edge node [midway, below=2pt]{\scriptsize \textbf{0~}} (yt11);
            \draw[thick, -latex](yt2) edge[bend left=60] node [midway, right=2pt]{\scriptsize \textbf{0~}}  (yt22);
            \draw[thick, -latex](ypt2) edge node [midway, below=2pt]{\scriptsize \textbf{0~}} (yt22);
            \draw[thick, -latex](yt2) edge[bend left=-5] node [midway, left=2pt]{\scriptsize \textbf{0~}}  (yt33);
            \draw[thick, -latex](ypt3) edge node [midway, below=2pt]{\scriptsize \textbf{0~}} (yt33);
        \end{tikzpicture}
        \caption{Transformed instance.}
    \end{subfigure}
    \caption{Transformation of QSTP. For each potential terminal $\vertexNumbered{i} \in \potTerminals$ (blue squares) a fixed terminal (red square) $\vertexNumbered{i^\prime}$ is added with profit $\vertexProfitsNumbered{i^\prime} = \vertexProfitsNumbered{i}$, which is connected to the root node $r$ and the original vertex by an arc of zero costs. The original potential terminal $\vertexNumbered{i}\in \potTerminals$ is removed from the set of potential terminals and considered a Steiner node (black circle) from here on. The original arc costs are omitted.}
    \label{fig:GraphTrans}
\end{figure}
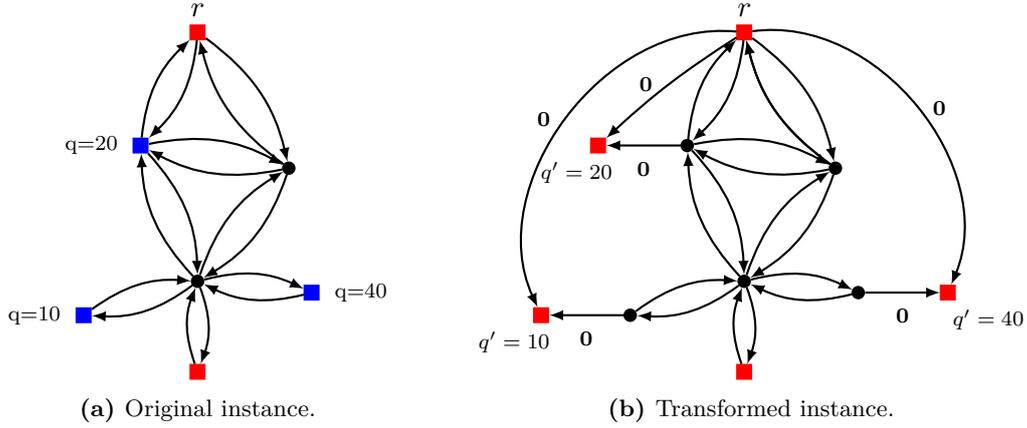

Now, let $\setArcs_r$ denote the set of arcs from the root node $r$ to the newly-added fixed terminals $\fixTerminalsTrans$ and $\setArcs_{t}$ the set of arcs from the original potential terminals \potTerminals to the newly-added fixed terminals $\fixTerminalsTrans$. Let $\setArcs^{\prime} = \setArcs \cup \setArcs_r \cup \setArcs_t$. Formulated as an IP, the problem reads as follows:
\begin{align}
    \min \quad&c^T x &\label{eq:obj_trans_qstp}\\
    \text{s.t.} \quad&&\nonumber\\
    & x(\delta^-(W)) \ge 1& \forall W \subset \setVertices^{\prime}, r\notin W, |W \cap \terminals^{\prime}| \ge 1\label{Steiner cut trans}\\
    & \sum_{\vertexNumbered{i^\prime} \in \fixTerminalsTrans} \vertexProfitsNumbered{i^{\prime}} x_{r, i^{\prime}} \le \sum_{\vertexNumbered{i^\prime} \in \fixTerminalsTrans} \vertexProfitsNumbered{i^{\prime}} - \quota & \label{eq:quota_cons_trans}\\
    & x_{ij} \in \lbrace 0,1 \rbrace &\forall (i,j) \in \setArcs^{\prime}\label{eq:var_trans_qstp}
\end{align}
Let \fullTransQSTPInstance denote an instance of the transformed QSTP. Figure \ref{fig:GraphTrans} shows the original graph and its transformation.

\subsubsection{Equivalence of the LP-relaxations}
We show the equivalence of the LP-relaxations of the \quotaInstance and the \quotaInstanceTrans, which are denoted by \lpRelaxation{\quotaInstance} and \lpRelaxation{\quotaInstanceTrans}, respectively. Let \optSolution{\generalProblem} denote the value of the optimal solution of a mathematical programming formulation \generalProblem and let $\mathcal{P}_{\lpRelaxation{}}(\generalProblem)$ denote the optimal solution of its LP relaxation. Before we compare the two LP-relaxations, let us introduce an additional set of variables in the transformed instance. For each vertex $\vertexNumbered{i}$ from the original set of potential terminals \potTerminals, let $\nodeVariable^\prime_i$ be defined as:
\begin{align}
        \nodeVariable_i^\prime \coloneqq\arcVariable(\incomingArcs{\vertexNumbered{i}}).
        \label{eq:nodevariable_extra}
\end{align}
Let this extended formulation be denoted as $\overline{\quotaInstanceTrans} = (\setVertices^{\prime}, \setArcs^{\prime}, \terminals^{\prime}, c^{\prime}, q^{\prime}, Q)$. By construction this new variable does not change the solution of the transformed instance $\quotaInstanceTrans$. However, we can now formulate the result of this section:

\begin{proposition}\label{prop:LP-equal}
Let \quotaInstance be an instance of the QSTP and $\overline{\quotaInstanceTrans}$ its transformed problem with the additional set of variables. It holds that:
\begin{align}
    proj_{xy}(\mathcal{P}_{\lpRelaxation{}}(\overline{\quotaInstanceTrans})) = \mathcal{P}_{\lpRelaxation{}}(\quotaInstance).
\end{align}
\end{proposition}

\begin{proof}
    \textbf{1) $\boldsymbol{proj_{xy}(\mathcal{P}_{\lpRelaxation{}}(\overline{\quotaInstanceTrans})) \subseteq\mathcal{P}_{\lpRelaxation{}}(\quotaInstance)}$:}
    let $(\arcVariable, \arcVariable^{t}, \arcVariable^{r}, y^\prime) \in \mathcal{P}_{\lpRelaxation{}}(\overline{\quotaInstanceTrans})$, where $\arcVariable$, $\arcVariable^{t}$, and $\arcVariable^{r}$ represent the variables for the arcs in $\setArcs$, $\setArcs_t$, and $\setArcs_r$, respectively. For readability, we drop the superscripts $t$ and $r$ when referring to a single variable $\arcVariable_{ii^\prime}$ of the arc $(i, i^\prime) \in A_t$ and $\arcVariable_{ri^\prime}$ of an arc $(r,i^\prime)\in A_r$, respectively. It is easy to see that $\arcVariable$ satisfies constraint \eqref{Steiner cut} as this must hold for every original terminal $t\in\fixTerminals$ in \lpRelaxation{\quotaInstanceTrans}, too. Now we show that $\nodeVariable^\prime$ satisfies constraints \eqref{Steiner cut pot} and \eqref{quota cons}.
    \begin{itemize}
        \item[(I)] $\nodeVariable^\prime$ satisfies \eqref{Steiner cut pot}: For each $\vertexNumbered{i^{\prime}} \in \fixTerminalsTrans$ consider the set $W^\prime \subset V^\prime$ with $r\notin W^\prime$ and $\lbrace\vertexNumbered{i}, \vertexNumbered{i^{\prime}} \rbrace \subseteq W^{\prime}$. We have
        \begin{align}
            \arcVariable(\incomingArcs{W^{\prime}}) = \arcVariable_{ri^{\prime}} + \arcVariable(\incomingArcs{W^{\prime}\backslash\vertexNumbered{i^{\prime}}}) \ge\arcVariable(\incomingArcs{W^{\prime}\backslash\vertexNumbered{i^{\prime}}}). \label{eq:proof_step1.1}
        \end{align}
        \textcite{goemans1993} showed that $\arcVariable(\incomingArcs{W^{\prime}\backslash\vertexNumbered{i^{\prime}}}) \ge \arcVariable(\incomingArcs{w})$ for any vertex $w \in W^\prime$, which is not a terminal, i.e., $w\notin T_f$. Thus, \ref{eq:proof_step1.1} gives
        \begin{align}
            \arcVariable(\incomingArcs{W^{\prime}\backslash\vertexNumbered{i^{\prime}}}) \ge \arcVariable(\incomingArcs{\vertexNumbered{i}}) \overset{\eqref{eq:nodevariable_extra}}{=} \nodeVariable_i^\prime \label{eq:proof_step1.2}
        \end{align}
        As by construction there solely exists an arc $(i,i^{\prime})$ and no arc $(i^{\prime},i)$ and \eqref{eq:proof_step1.2} holds for every $W^\prime\backslash\vertexNumbered{i^{\prime}}$, it also holds for the set $W$ in the original graph with $W\subset \setVertices$, $r\notin W$, and $\vertexNumbered{i}\in W$. Thus, we have
        \begin{align}
            &\arcVariable(\incomingArcs{W}) \ge \nodeVariable_i^\prime & \forall W \subset \setVertices, r\notin W, |W \cap \potTerminals| \ge 1, \vertexNumbered{i} \in \potTerminals, \label{eq:proof_step1.3}
        \end{align}
        which equals \eqref{Steiner cut pot}. 
        \item[(II)] $\nodeVariable^\prime$ satisfies \eqref{quota cons}: Consider the set $W^\prime = \lbrace\vertexNumbered{i}, \vertexNumbered{i^{\prime}} \rbrace$ for each $\vertexNumbered{i^{\prime}} \in \fixTerminalsTrans$. As $W^\prime$ contains a fixed terminal and does not contain the root, the set has to have at least one incoming arc and \eqref{Steiner cut trans} applies:
        \begin{align}
            \arcVariable(\incomingArcs{W^{\prime}} = \arcVariable_{ri^{\prime}} + \arcVariable(\incomingArcs{\vertexNumbered{i}}) \ge 1 &&\Leftrightarrow &&\arcVariable_{ri^{\prime}} \ge 1 - \arcVariable(\incomingArcs{\vertexNumbered{i}}). \label{eq:proof_step2.1}
        \end{align}
        Given the quota constraint \eqref{eq:quota_cons_trans} of the transformed QSTP:
        \begin{align*}
        \sum_{\vertexNumbered{i^{\prime}} \in \fixTerminalsTrans} \vertexProfitsNumbered{i^{\prime}} - \quota \ge \sum_{\vertexNumbered{i^{\prime}} \in \fixTerminalsTrans} \vertexProfitsNumbered{i^{\prime}} \arcVariable_{ri^{\prime}} \overset{\eqref{eq:proof_step2.1}}{\ge} \sum_{\vertexNumbered{i^{\prime}} \in \fixTerminalsTrans} \vertexProfitsNumbered{i^{\prime}} (1 - \arcVariable(\incomingArcs{\vertexNumbered{i}})) \\
        \Leftrightarrow -\quota \ge - \sum_{\vertexNumbered{i^{\prime}} \in \fixTerminalsTrans} \vertexProfitsNumbered{i^{\prime}} \arcVariable(\incomingArcs{\vertexNumbered{i}}) \\
        \Leftrightarrow \quota \le \sum_{\vertexNumbered{i^{\prime}} \in \fixTerminalsTrans} \vertexProfitsNumbered{i^{\prime}} \arcVariable(\incomingArcs{\vertexNumbered{i}}) \overset{\eqref{eq:nodevariable_extra}}{=} \sum_{\vertexNumbered{i} \in \potTerminals} \vertexProfitsNumbered{i}\nodeVariable^\prime_i .
        \end{align*}
        Hence, $y^\prime$ satisfy \eqref{quota cons}. 
    \end{itemize}
    Finally, by ignoring $\arcVariable^r$ and $\arcVariable^t$, we have $proj_{xy}(\mathcal{P}_{\lpRelaxation{}}(\overline{\quotaInstanceTrans}))~\subseteq~\mathcal{P}_{\lpRelaxation{}}(\quotaInstance)$. This concludes the first part of the proof.

    \textbf{2) $\boldsymbol{proj_{xy}(\mathcal{P}_{\lpRelaxation{}}(\overline{\quotaInstanceTrans})) \supseteq\mathcal{P}_{\lpRelaxation{}}(\quotaInstance)}$:} Given a feasible solution $(\hat{\arcVariable}, \hat{\nodeVariable}) \in \mathcal{P}_{\lpRelaxation{}}(\quotaInstance)$, we construct $(\arcVariable, \arcVariable^{t}, \arcVariable^{r}, y^{\prime}) \in \mathcal{P}_{\lpRelaxation{}}(\overline{\quotaInstanceTrans})$. For all $(i,j) \in A$, set $\arcVariable_{ij} = \hat{\arcVariable}_{ij}$. For all $\vertexNumbered{i} \in \potTerminals$, set $\arcVariable_{ii^{\prime}} = y_i^{\prime} = \hat{\nodeVariable}_i$ for all $(i, i^\prime) \in \setArcs_t$ and $\arcVariable_{ri^{\prime}} = 1 - \hat{\nodeVariable}_i$ for all $(r,i^{\prime})\in\setArcs_r$. In the following we show that $(\arcVariable, \arcVariable^{t}, \arcVariable^{r})$ satisfy the Steiner cut constraint \eqref{Steiner cut trans} and the quota constraint \eqref{eq:quota_cons_trans} of the transformed instance: 
    \begin{itemize}
        \item[(I)] For any original fixed terminal $\arcVariable$ satisfies \eqref{Steiner cut trans}. For any newly-added fixed Terminal $\vertexNumbered{i^{\prime}} \in \fixTerminalsTrans$, consider the Steiner cut constraint \eqref{Steiner cut trans} of the set $W^\prime=\lbrace\vertexNumbered{i^\prime}\rbrace$:
    \begin{align*}
        \arcVariable(\incomingArcs{\vertexNumbered{i^{\prime}}}) = \arcVariable_{ri^{\prime}} + \arcVariable_{ii^{\prime}} = 1 - \hat{\nodeVariable}_i + \hat{\nodeVariable}_i = 1 \ge 1
    \end{align*}
    Now, for every $W^\prime \subset \setVertices$ with $r\notin W^\prime$, $\lbrace \vertexNumbered{i}, \vertexNumbered{i^\prime}\rbrace \in W^\prime$ and $W = W^{\prime}\backslash\vertexNumbered{i^{\prime}}$, the constraint \eqref{Steiner cut trans} yields:
    \begin{align*}
        \arcVariable(\incomingArcs{W^\prime}) =  \arcVariable_{ri^{\prime}} + \arcVariable(\incomingArcs{W^\prime\setminus\vertexNumbered{i^\prime}}) = 1 - \hat{\nodeVariable}_i + \arcVariable(\incomingArcs{W}) \overset{\eqref{Steiner cut pot}}\ge 1 - \hat{\nodeVariable}_i + \hat{\nodeVariable}_i \ge 1
    \end{align*}
    So, $(\arcVariable, \arcVariable^{t}, \arcVariable^{r})$ satisfy \eqref{Steiner cut trans}.
    
    \item[(II)] As $(\hat{\arcVariable}, \hat{\nodeVariable}) \in \mathcal{P}_{\lpRelaxation{}}(\quotaInstance)$, the quota constraint \eqref{quota cons} is satisfied:
    \begin{align*}
        \quota \le \sum_{\vertexNumbered{i} \in \potTerminals} \vertexProfitsNumbered{i}\hat{\nodeVariable}_i =  \sum_{\vertexNumbered{i^{\prime}} \in \fixTerminalsTrans} \vertexProfitsNumbered{i^\prime} (1 - \arcVariable_{ri^{\prime}}) =  \sum_{\vertexNumbered{i^{\prime}} \in \fixTerminalsTrans} \vertexProfitsNumbered{i^\prime} - \sum_{\vertexNumbered{i^{\prime}} \in \fixTerminalsTrans} \vertexProfitsNumbered{i}\arcVariable_{ri^{\prime}} \\
        \Leftrightarrow \sum_{\vertexNumbered{i^{\prime}} \in \fixTerminalsTrans}  \vertexProfitsNumbered{i^\prime}\arcVariable_{ri^{\prime}} \le \sum_{\vertexNumbered{i} \in \potTerminals} \vertexProfitsNumbered{i^\prime} - \quota
    \end{align*}
    So, $(\arcVariable, \arcVariable^{t}, \arcVariable^{r})$ satisfy \eqref{eq:quota_cons_trans}. 
    \end{itemize}
    With (I) and (II) it is shown that $(\arcVariable, \arcVariable^{t}, \arcVariable^{r}, y^{\prime}) \in\mathcal{P}_{\lpRelaxation{}}(\quotaInstanceTrans)$. This concludes the second part of the proof.
   
   With $proj_{xy}(\mathcal{P}_{\lpRelaxation{}}(\overline{\quotaInstanceTrans}))~\subseteq~\mathcal{P}_{\lpRelaxation{}}(\quotaInstance)$ and $proj_{xy}(\mathcal{P}_{\lpRelaxation{}}(\overline{\quotaInstanceTrans}))~\supseteq~\mathcal{P}_{\lpRelaxation{}}(\quotaInstance)$, we have $proj_{xy}(\mathcal{P}_{\lpRelaxation{}}(\overline{\quotaInstanceTrans})) = \mathcal{P}_{\lpRelaxation{}}(\quotaInstance)$ which concludes the proof of Proposition \ref{prop:LP-equal}
\end{proof}

\subsection{Bi-objective QSTP}\label{sec:BiObj}
As mentioned in the introduction, not only costs but also other social or environmental impacts must be taken into account when planning new wind farms. In this study, we focus on minimizing both costs and the impact on the landscape of network cables and wind turbines. In addition to the real cable (edge) costs $\cost$ and wind turbine costs $\vertexCosts$ introduced in Section \ref{sec:QSTP}, let $\scenicness(\edge)$ and $\scenicnessVertex(\vertex)$ denote the scenic impact of an edge $\edge \in \setEdges$ and of a vertex $\vertex\in \potTerminals$, respectively. The objective is to minimize both the costs $C(S)$ and the scenic impact $L(S)$ of the Steiner tree $S = (\setEdges ^\prime, \setVertices^\prime) \subseteq \graph$, i.e.,
\begin{align}
    C(S) = \sum_{\edge\in \setEdges^\prime} \cost(\edge) + \sum_{\vertex \in \potTerminals \cap \setVertices^\prime} \vertexCosts(\vertex) &\quad \mathrm{and} & L(S) = \sum_{\edge\in \setEdges^\prime} \scenicness(\edge) + \sum_{\vertex \in \potTerminals \cap \setVertices^\prime} \scenicnessVertex(\vertex)
\end{align}
while fulfilling the quota $\quota$, see \eqref{eq:quota_definition}.

We are interested in finding Pareto optimal solutions to this problem. A solution is Pareto optimal (non-dominated), if there exists no other feasible solution with a lower objective value for both goals \parencite{ehrgott2016}. Usually, finding all Pareto optimal solutions is computationally hard. Multiple ways to approach such a multi-objective optimization problem exist, each of them having its own advantages and disadvantages. Among others, the common methods are, e.g., the weighted-sum approach and the $\varepsilon$-constraint method; see 
\parencite{ehrgott2016} for a general overview. In the context of STP, \textcite{leitner2015} proposes an $\varepsilon$-constraint algorithm to solve a bi-objective PCSTP. Although the $\varepsilon$-constraint method produces a large number of solutions, many of them are dominated by others. There have been efforts to find representative subsets of non-dominated solutions, see, e.g., \textcite{dogan2022, ceyhan2019, mesquita-cunha2023}. However, since the single variant of the integrated wind layout and cable routing problem, discussed in this study, is already hard to solve, we restrict ourselves to the weighted-sum approach. Despite only generating supported solutions\footnote{supported solutions: points on the convex relaxation of the Pareto frontier}, its advantage is that, in contrast to other methods, such as the $\varepsilon$-constraint one, we can retain the general structure of the single-objective QSTP formulation \parencite{ehrgott2016}. In this study, we compose a new objective function $\overline{C}$ by using a convex combination of the real costs $\cost$ and $\vertexCosts$ and the scenic impact $\scenicness$ and $\scenicnessVertex$ as follows:
\begin{equation}
    \overline{C}(S) = \alpha C(S) + (1 - \alpha) L(S),
\end{equation}
where $\alpha \in [0,1]$. Depending on the choice of $\alpha$, a different solution is found, with $\alpha=1.0$ only taking costs into account and $\alpha=0.0$ only considering the scenic impact.

\subsection{Shortest path reduction and heuristic}

In general, reduction techniques reduce the size of the original problem by removing arcs and vertices without cutting the optimal solution. Reduction techniques are often used in a preprocessing step to create "easier to solve" instances.

A common and intuitive reduction test is the shortest path reduction: If there exists a directed path $P(v,w) = \lbrace v, (v,v_1), v_1, (v_1,v_2), \dots, v_n, (v_n, w), w\rbrace$ of costs $c(P(v,w))$ with
\begin{align*}
    c(P(v,w)) < c(v,w)
\end{align*}
then arc $(v,w)$ can be removed. 

\begin{proposition}
Consider an instance \fullTransQSTPInstance of the previously-described QSTP. If there exists a directed path $P(v,w)$ with costs $c(P(v,w)) < c(v,w)$, then there also exists a directed path $P(w,v)$ with costs $c(P(w,v)) < c(w,v)$, and both $(v,w)$ and $(w,v)$ can be removed from $I_{QT}$.
\end{proposition}

\begin{proof}
Given the path $P(v,w)$ with
\begin{align*}
    c(P(v,w)) < c(v,w)
\end{align*}
each node $v_i \in \lbrace v, v_1, \dots, v_n, w\rbrace$ on the path $P(v,w)$ is either in $T_p$, i.e, $w_{v_i} > 0$, or not, i.e. $w_{v_i} = 0$. Now, subtracting the costs of node $w$, $w_w$, and adding the costs of node $v$, $w_v$, yields:
\begin{align}
    c(P(v,w)) - w_w + w_v < c(v,w) - w_w + w_v \nonumber\\
    c(v,v_1) + c(v_1,v_2) + \ldots + c(v_n, w) - w_w + w_v < c(v,w)- w_w + w_v \label{shortest path costs}
\end{align}
The cost of an arc $(i,j)$ is given by $c(i,j) = c_e + w_j$ where $c_e$ is the costs of the undirected edge $e=(i,j)$. With this, \eqref{shortest path costs} reads as follows:
\begin{align*}
    &\overbrace{c_{v,v_1} + w_{v_1}}^{c(v,v_1)} + \overbrace{c_{v_1,v_2} + w_{v_2}}^{c(v_1, v_2)} + \ldots + \overbrace{c_{v_n, w} + w_w}^{c(v_n,w)} - w_w + w_v < c_{v,w} + w_w - w_w + w_v \\
    &\underbrace{w_v + c_{v,v_1}}_{c(v_1, v)} + \underbrace{w_{v_1} + c_{v_1,v_2}}_{c(v_2,v_1)} + \ldots + \underbrace{w_{v_n} + c_{v_{n}, w}}_{c(w,v_n)} < c_{v,w} + w_v \\
    &c(P(w,v)) < c(w,v)
\end{align*}
Hence, there also exists a path $P(w,v)$ with less costs as the costs of arc $(w,v)$.
\end{proof}

Heuristics are used to find good and feasible solutions in short time, and thus providing upper bounds on the exact solution of a problem. In the context of STP, one of the best-known heuristics is the \textit{shortest path heuristic}, which was introduced by \textcite{takahashi1980}. For this study, we adapt the implementation of the shortest path heuristic given by \textcite{rehfeldt2021a}. Starting at the root node, we add the closest fixed terminal or potential terminal to the root component. Then we update the distance of all non-connected vertices to the root component. We repeat until all fixed terminals are connected and the connected potential terminals fulfill the desired quota.

\section{Computational results}\label{sec:Results}

In this section, we describe the computational implementation and setup in \ref{sec:implement}, the input data in \ref{sec:input_sets}, and the computational results in \ref{sec:compRes} to \ref{sec:Res_Sequentiel_vs_combined}.

\subsection{Implementation details}\label{sec:implement}

In the following, we roughly describe the functionality of the exact Steiner tree solver \scipjack. We present a flow-based MIP formulation which we use to benchmark our proposed method. Subsequently, we introduce the setting used for all computations. 

\subsubsection{Extension of SCIP-Jack}

The transformed QSTP (\eqref{eq:obj_trans_qstp} -- \eqref{eq:var_trans_qstp}) is integrated in the general STP solver \scipjack. \scipjack uses a branch-and-cut procedure to account for the exponential number of constraints induced by the Steiner cut-like constraints, i.e., \eqref{Steiner cut trans}. We use the native, flow-based separation algorithm of \scipjack for the Steiner-cut generation and add the additional quota constraint. Most other features of \scipjack, such as heuristics, reduction methods, and domain propagation, do not work with the quota constraint and must be omitted. Therefore, we implemented the shortest-path-based reduction method, which is used in the presolving step. We implement the shortest-path heuristic as primal heuristic to improve the solving process. The heuristic might find an optimal solution or gives a primal bound, with which nodes of the branch-and-bound tree can be removed. Depending on the current LP-solution $\hat{\arcVariable}$, the shortest-path heuristic is called repeatedly on a graph with modified arc costs $\overline{\cost}_\arc = (1.0 - \hat{\arcVariable}_\arc) \cost_\arc$ for all $\arc \in \setArcs$ during the branch-and-bound procedure.

\subsubsection{Flow-based MIP formulation}\label{sec:flow-based}
We verify our proposed formulation by solving the problem with the general out-of-the-box MIP solver \gurobi \parencite{gurobi95}. Due to the exponential number of constraints given by \eqref{Steiner cut trans}, which cannot be handle by general MIP solvers, we model the problem as the following flow-based MIP formulation (\texttt{FLOW}). Let $r \in \fixTerminals$ be chosen as root. Note that in our case as all substations (terminals) are interconnected by edges of zero costs, it does not matter which terminal is chosen as the root. The flow-based MIP formulation is given by
\begin{align}
    \min \quad&\cost^T \arcVariable + \vertexCosts^T y& \label{eq:obj_flow}\\
    \text{s.t.} \quad&&\\
    &\sum_{\vertex\in\potTerminals} \vertexProfitsNumbered{\vertex} y_\vertex \ge \quota \label{eq:quota_flow}\\
    &\sum_{\arc\in\incomingArcs{\vertex}}\arcFlow_\arc - \sum_{\arc\in\outgoingArcs{\vertex}}\arcFlow_\arc = \begin{dcases} 0 & \forall \vertex \in \setVertices\setminus(\fixTerminals \cup\potTerminals)\\ 
    1 & \forall \vertex \in \fixTerminals\setminus r\\
    y_\vertex & \forall \vertex \in \potTerminals
    \end{dcases} \label{eq:flowBalance_flow}\\
    &\arcVariable_\arc \le y_\vertex &\forall \arc \in \incomingArcs{\vertex}, \forall \vertex \in \potTerminals \label{eq:activeTerm_flow}\\
    &\arcFlow_\arc \le M \arcVariable_\arc &\forall \arc \in \setArcs \label{eq:activeArc_flow}\\
    &y_\vertex \in \lbrace 0,1\rbrace &\forall\vertex \in \potTerminals\\
    &\arcVariable_\arc \in \lbrace 0,1\rbrace, \,\arcFlow_\arc \in \R_{\ge0} &\forall\arc \in \setArcs \label{eq:var_flow}
\end{align}

where $\arcVariable$ and $y$ denote the decision variables if an arc and a vertex is chosen, respectively, and $\arcFlow$ describes the flow over the arcs. Constraint \eqref{eq:quota_flow} describes the quota constraint. Constraint \eqref{eq:flowBalance_flow} captures the flow balance at each vertex depending on its type. For any fixed terminal, in our case the grid's substations have to be connected to the root node, so the right-hand side is equal to one, so a "flow" $\arcFlow$ of one has to reach the terminal. Also, any chosen potential terminal (turbine), i.e., $y_\vertex = 1$, has to be connected to the root node. The incoming arcs of a potential terminal can only be active if the potential terminal is chosen, as in \eqref{eq:activeTerm_flow}. Equation \eqref{eq:activeArc_flow} ensures that a flow over an arc is only possible if the arc is active. The big-M notation is used to limit the flow on an active arc, i.e., the upper bound on $\arcFlow_\arc$. We choose $M = | \fixTerminals \cup \potTerminals |$ as an upper limit, which would allow all fixed and potential terminals to be connected via a single string. In terms of regional wind farm planning, this limit is an overestimation, as the potential wind turbines are usually distributed around multiple substations, which are already connected via the existing power grid. 

\subsubsection{Computational setting}\label{sec:comp_setting}
The flow-based formulation was implemented in \python~\oldstylenums{3.8.10} using the \gurobi~\python-interface, and is solved with \gurobiVersion{9}{5}~\parencite{gurobi95}. For the QSTP and its transformation we use \scipjack in \scipVersion{8}{0}{1} \parencite{bestuzheva2023a} using \cplexVersion{12}{10} \parencite{ibmilog2022} as the LP solver.

We implemented the following settings: a) the flow-based formulation (\texttt{FLOW} with \gurobi), b) the initial QSTP formulation \texttt{QSTP}, c) the transformed QSTP \texttt{TransQSTP}, d) the transformed QSTP plus a shortest-path reduction \texttt{TransQSTP+}, and e) the transformed QSTP plus the shortest-path reduction and a shortest-path based heuristic \texttt{TransQSTP++}. All computations were executed single-threaded in the case of \scipjack and using up to 32 threads in the case of \gurobi in non-exclusive mode on a cluster with \textit{Intel XeonGold 6342} CPUs running at 2.8 GHz, where five CPUs and 32 GB of RAM is reserved for each run. We set a time limit of six hours (21,600 s).

\subsection{Input data and instance sets}\label{sec:input_sets}

The considered wind turbines are potential turbines with costs and designs for the year 2050 from \textcite{Ryberg.2019}.
To determine their locations, state-of-the-art methods were applied to exclude unsuitable areas, for example due to high terrain steepness or minimal distances to settlements or infrastructure. A total of around 160,000 turbine sites were then selected in the remaining suitable areas in Germany, taking into account minimum distances between turbines based on wind roses, etc. This means that the potential turbine sites from which a certain number of turbines could be selected are fixed in the optimization. \textcite{Ryberg.2019} obtained the costs for the installed wind turbines based on future turbine designs for 2050 as well as their annual energy yields. In calculating the energy yields, some stochastic effects that contribute to losses such as turbine availability, electrical inefficiencies or wake effects have already been covered by the convolution of the power curve. In a validation step with measured data in \textcite{Ryberg.2019}, over-prediction of power in low-generation times has been addressed by implementing a capacity-factor-dependent correction factor. The specific costs for cables and the grid connection of turbines are based on \textcite{McKenna.2021b}. The data for substation locations is based on OpenStreetMap \parencite{OpenStreetMapcontributors.2021} entries. 

\begin{figure}
	\centering
	\includegraphics[width=0.6\columnwidth]{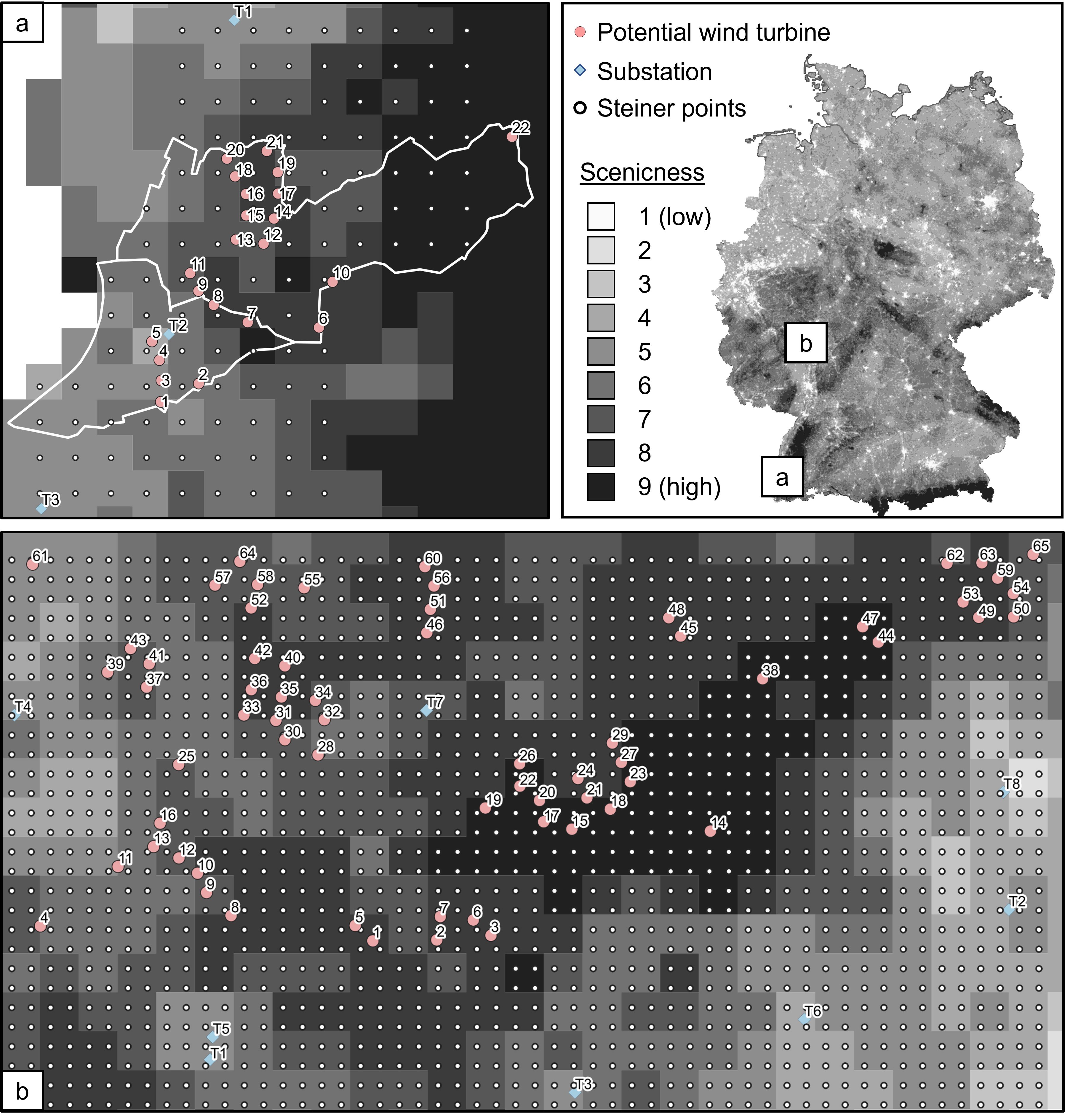}
	\caption{Substations, potential wind turbines, Steiner points, and scenic value in case studies A and B in Germany. The map of Germany shows the locations of the two case studies as well as the landscape's scenic value distribution in Germany.}
	\label{fig:overview problem}
\end{figure}

The scenic features of landscapes in Germany are derived from \textcite{Roth.2018}, 
who for the first time conducted an area-wide analysis of the aesthetic value of landscapes across the country. Using selected landscape photographs, a survey and subsequent regression, each 1 km$^2$ in Germany was assigned a discrete aesthetic value from one (low scenic quality) to nine (high scenic quality) (see Figure \ref{fig:overview problem}). \textcite{Roth.2018} ensured the representativeness of the landscape images for the entire federal territory by selecting 30 reference spaces with an area of 130 to 140 km$^2$. The photographs were then tagged with metadata such as whether they show forests, water bodies or settlements. In the subsequent representative survey, the participants were asked to rate the scenic beauty of landscapes in 10 randomly selected photographs from a pool of 822 images. In a final step, the landscape values were extrapolated to the entire German territory based on the metadata using a regression analysis with 17 significant variables and a coefficient of determination of 0.639. For the showcase of our methodology, we assume that one km of electricity network has the same impact on the landscape as one wind turbine.

\begin{figure}
    \centering
    \includegraphics[width=.55\columnwidth]{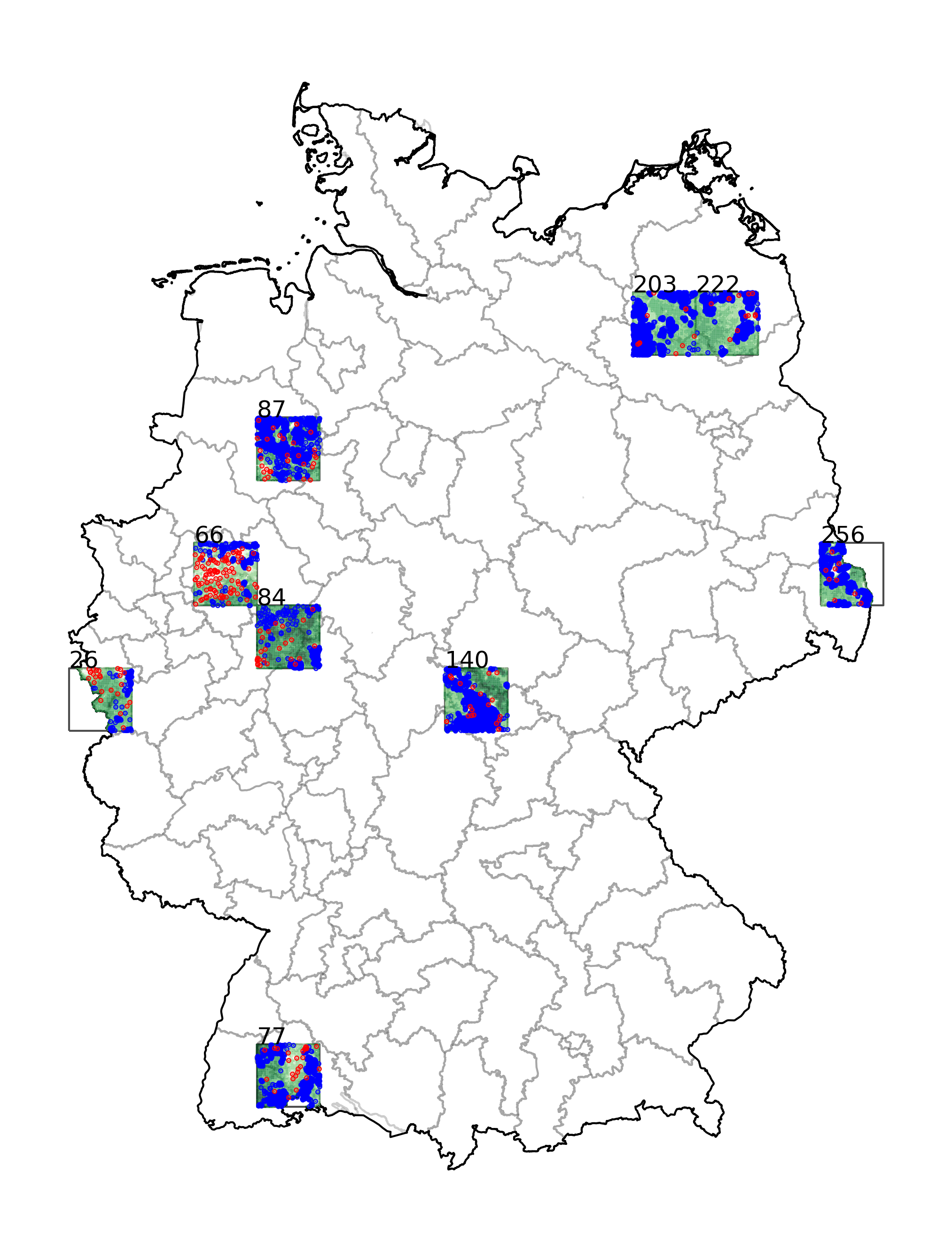}
    \caption{Chosen cells (numbered) overview; Turbines: blue; Substations: red}
    \label{fig:cells_overview}
\end{figure}

\begin{figure}
    \centering
    \includegraphics[width=.7\columnwidth]{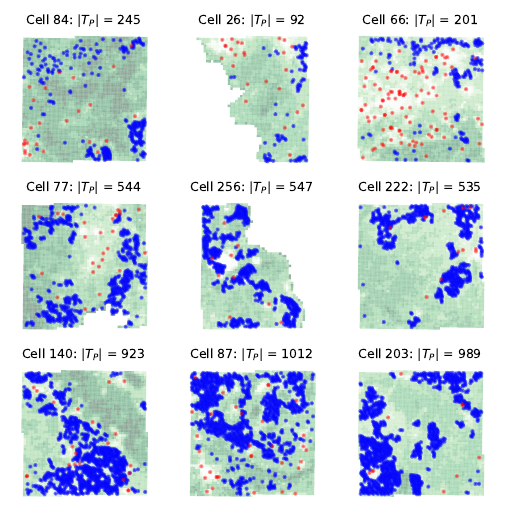}
    \caption{Chosen cells with number of turbines $|\potTerminals|$; Turbines: blue; Substations: red}
    \label{fig:cells_detailed}
\end{figure}

To evaluate our proposed approach we generate two benchmark sets using the data described above. In the first benchmark set a general choice of regions is made focusing on the number of turbines. For the second, a particular choice of two municipalities is made including a high number of Steiner nodes and, thus, a more flexible cable routing. In these municipalities the trade-offs between costs and scenic impact are evaluated in \ref{sec:tradeoffs}. For the first one, we evenly divide Germany into areas of 50~km x 50~km. We define three sets of sizes based on the number of potential turbine locations, i.e., \texttt{small} 50 to 300, \texttt{medium} 450 to 550, \texttt{large} 900 to 1100 turbines. For each of the categories, we choose the areas with minimal, median, and maximal scenic features of landscape, see Figure \ref{fig:cells_overview} and \ref{fig:cells_detailed}. For each of these areas we consider two topologies: 1) The complete graph $\graphFull$, where \setVertices consists of the set of substations \fixTerminals and the set of potential wind turbines \potTerminals allowing the cable routing only through chosen turbines; 2) Each turbine location is connected to an artificial Steiner node $v$ by an edge of zero costs. The set of edges consists of these edges and the edges of the complete subgraph $\graph^\prime = (\setVertices\setminus\potTerminals, E^\prime)$. Thus, each turbine position can be used for the cable routing without choosing to build the turbine. This results in six topologies over all sizes, three without and three with additional Steiner nodes. For each of the topologies, we create instances with five values of quota \quota, namely, 25\%, 40\%, 50\%, 60\%, 75\% of the maximal potential of the area, and using seven values for the weighting factor $\alpha$, namely 0.0, 0.2, 0.4, 0.5, 0.6, 0.8, 1.0, with $\alpha = 1.0$ considering only costs and $\alpha=0.0$ considering only the scenic value, to create the edge and vertex costs as described in Section \ref{sec:BiObj}.

For the second benchmark set, we take a more detailed look at two regions in Germany (see Figure \ref{fig:overview problem}). The first of these, region A, includes the siting and grid connection planning of potential wind turbines in the municipalities of Bad Bellingen and Schliengen. In these two municipalities, a total of 22 potential turbines can be connected to three different substations with a total annual energy yield of 158 GWh. For a flexible cable routing, 122 Steiner points are placed on a grid of dimension 1000 meters. A choice can be made between 10,731 possible edges. For region B, which contains 65 possible turbines with a total energy yield of 405 GWh and eight substations, a large instance is created by adding 1566 Steiner points on a grid with a dimension of 500 meters. A choice can be made between around 1.34 million edges.

For each of the two case study regions, we generate a set of instances on which we evaluate our proposed techniques. We vary both the given quota \quota and weight of the cable costs and scenic value $\alpha$, resulting in different edge costs. For region A, we choose the following quotas $\quota \in [10,20,30,40,50,60,70,80,90,100,110,120,130,140,150,158]$ and for region B we choose $\quota \in [50,100,150,175,200,250,300,350,405]$. For both regions, we let $\alpha$ increase from 0.0 to 1.0 in 0.1 steps to create the edge costs using the weighted-sum approach (Section \ref{sec:BiObj}). This results in 176 instances for region A and 99 for region B. Combining both of these benchmark sets yields 905 instances with different number of turbines, substations and density of Steiner nodes for the cable routing\footnote{The raw input data and the resulting QSTP instance files see \url{https://doi.org/10.5281/zenodo.10600575}}. A summary of the testsets is given in Table \ref{tab:testset_description}. 

\begin{table}[]
    \footnotesize
    \centering
    \caption{\normalsize Summary of instances. The third column gives the number of turbines (potential terminals) $|\potTerminals|$, the fourth column shows the range of number of nodes $|\setVertices|$ consisting of turbines, substations (fixed terminals), and possible Steiner nodes, and the fifth column shows the range of number of edges that can be chosen for the cable routing.}
    \begin{tabular}{lccccc>{\centering\arraybackslash}p{.01\linewidth}}
    Name & Instances & $|\potTerminals|$ & $|\setVertices|$ & $|\setEdges|$ & Description\\
    \hline\\
    Small  & 210 & 92, 201, 245& 117 - 515 & 6786 - 46056& \multirow{3}{*}{\begin{minipage}{0.25\linewidth}
        Areas of 50 km x 50 km chosen by number of potential turbines and value of scenicness, for each area: 5x $\quota$, 7x $\alpha$, without Steiner nodes and with Steiner nodes for every turbine
    \end{minipage}}\\[19pt]
    Medium & 210 & 535, 544, 547 & 547 - 1114 & 149331 - 162709 &\\[19pt]
    Large  & 210 & 923, 989, 1012 & 944 - 2057& 445096 - 546502 &\\[19pt]
    Region A    & 176 & 22 & 147  & 10731   &\begin{minipage}{0.25\linewidth}Steiner nodes on a grid of 1 km, 16x~$\quota$, 11x $\alpha$\end{minipage}\\[10pt]
    Region B    & 99  & 65 & 1639 & 1342341 &\begin{minipage}{0.25\linewidth}Steiner nodes on a grid of 0.5 km, 9x~$\quota$, 11x $\alpha$\end{minipage}\\[10pt]
    \hline
    \end{tabular}
    \label{tab:testset_description}
\end{table}

\subsection{Comparison with state-of-the-art generic solvers}\label{sec:compRes}
\begin{figure}
    \includegraphics[scale=1.]
    {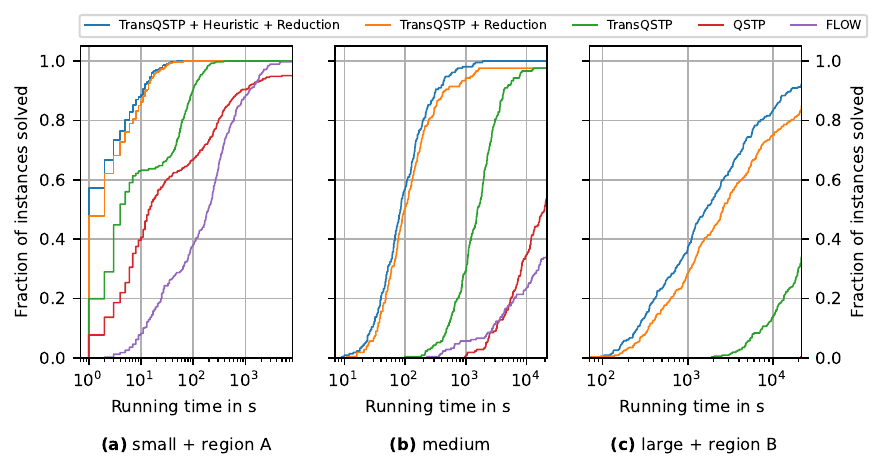}
    \caption{Performance profiles for \texttt{small}~+~\texttt{region A} (Panel (a)), \texttt{medium} (Panel(b)), and \texttt{large}~+~\texttt{region B} (Panel (c)), showing the cumulative percentage of instances solved to optimality over time for all approaches. The Running time in (a) and (b) is shown in seconds on a log-scale.}
    \label{fig:perfProfile_all_testsets}
\end{figure}

The summary of all computational results is presented in Table \ref{tab:testset_results_all}. Due to similar number of nodes and edges, region A and region B are classified as \texttt{small} and \texttt{large}, respectively. Figure \ref{fig:perfProfile_all_testsets} shows the fraction of instances that are solved to optimality over time for each testset class, i.e., (a) \texttt{small} + region A, (b) \texttt{medium}, and (c) \texttt{large} + region B.

In the case of \texttt{small} and region A, all 386 instances are solved to optimality by all settings but not the \texttt{QSTP}, which only solves 368 instances. The naive flow formulation solved by standard out-of-the-box MIP solver is significantly outperformed by an order of magnitude by \texttt{TransQSTP} and by two orders of magnitude by \texttt{TransQSTP++} and \texttt{TransQSTP+}, see Figure \ref{fig:perfProfile_all_testsets}(a). \texttt{QSTP} is also clearly dominated by its transformed version, even without any reduction or heuristic. Additionally, \texttt{TransQSTP++} and \texttt{TransQSTP+} instantly solve 57\% and 48\% of the instances, respectively, and almost 90\% in ten seconds. Even though \texttt{QSTP} solves around 40\% of the instances after 10s, it fails to solve all instances. With \texttt{FLOW}, it takes 8000s to solve all the instances to optimality. 

Considering the \texttt{medium}-sized instances in Figure \ref{fig:perfProfile_all_testsets}(b), the advantage of the specialized QSTP approach becomes even more evident. Whereas \texttt{TransQSTP++} and \texttt{TransQSTP+} show a similar performance with \texttt{TransQSTP+} leaving a very small gap ($<0.05\%$) on five instances, only half and a third of the instances are solved to optimality by \texttt{QSTP} and \texttt{FLOW}, respectively. Again, \texttt{TransQSTP++} and \texttt{TransQSTP+} outperfrom \texttt{FLOW} by two orders of magnitude and \texttt{TransQSTP} by an order of magnitude. Out of the instances that are not solved to optimality by \texttt{FLOW}, 53 instances are stopped by the time limit and 83 instances due to memory limit. 

The problem-specific reduction method and the sophisticated, problem-related, cut-seperation algorithm implemented in \scipjack are the main reasons for the advantages compared to the standard out-of-the-box solver. First, \scipjack can often drastically reduce the problem size by utilizing the structure of the problem graph itself, whereas the generic presolving methods of MIP solvers are usually remarkably unsuccesfull on Steiner tree instances. Second, \scipjack utilizes strong, but highly intricate, problem-specific cutting planes, whereas generic solvers are not able to produce such cutting planes, because they depend on a combinatorial knowledge of the underlying problem. Finally, another advantage, which will be shown next, is that \scipjack uses problem specific primal heuristics, whereas MIP solvers often struggle to find a primal solution of good quality.


So far \texttt{TransQSTP++} has only shown a small improvement with respect to \texttt{TransQSTP+}. It is common that the main improvement on solving Steiner tree related problems comes from presolving techniques, see e.g., \cite{rehfeldt2021a}. However, when looking at the results of the \texttt{large}-sized instances in Figure \ref{fig:perfProfile_all_testsets}(c), \texttt{TransQSTP++} finds more solutions and solves more of them to optimality compared to all other settings. More accurately, \texttt{TransQSTP++} solves 285 out of 309 instances to optimality within the time limit with an average gap of 4.8\% for the remaining ones while \texttt{TransQSTP+} solves only 262 instances to optimality and has larger gaps on the non-optimal instances it finds a solution for. \texttt{FLOW} does not find any optimal solution and only a solution on 290 instances with an average gap of 68.5\% (maximum gap > 1000\%). 128 instances run into memory limit and 181 instances hit the time limit.

\begin{table}[]
    \small
    \caption{\normalsize Summary of results. The table shows for each instance set and setting the number of instances for which a solution was found and that were solved to optimality in the fourth and fifth column, respectively. The last three columns show the average gap on the instances which were not solved to optimality, the average time needed to solve the instances to optimality, and the number on how many instances the setting was the fastest/among the fastest (best/wins).}
    \begin{tabular}{>{\arraybackslash}p{.1\linewidth}lcccrrr}
    Testset& Setting &	\# & \# Sol. & \# Opt.  &   $\varnothing$ Gap [\%]   & $\varnothing$ Time [s] & best/wins \\[0.5ex]
    \hline\hline
    \rule{0pt}{2.6ex}
    \multirow{5}{*}{\begin{minipage}{.9\linewidth}small + region A\end{minipage}}	 & TransQSTP++  &	 386 &	 386 &	 386 &	     --  &	     3.96 &	 127/333	 \\ 	 
			 & TransQSTP+   &	 386 &	 386 &	 386 &	     --  &	     4.70 &	 53/259	 \\ 	 
			 & TransQSTP    &	 386 &	 386 &	 386 &	     --  &	    32.63 &	 0/77	 \\ 	 
			 & QSTP         &	 386 &	 386 &	 368 &	    0.46 &	   206.83 &	 0/30	 \\ 	 
			 & FLOW         &	 386 &	 386 &	 386 &	     --  &	   440.06 &	 0/0	 \\[0.5ex] 	 
    \hline
    \rule{0pt}{2.6ex}
    \multirow{5}{*}{\begin{minipage}{.9\linewidth}medium\end{minipage}}	 & TransQSTP++  &	 210 &	 210 &	 210 &	     --  &	   157.83 &	 137/141	 \\ 	 
			 & TransQSTP+   &	 210 &	 210 &	 205 &	    0.05 &	   191.92 &	 69/73	 \\ 	 
			 & TransQSTP    &	 210 &	 210 &	 205 &	    0.04 &	  2043.05 &	 0/0	 \\ 	 
			 & QSTP         &	 210 &	 210 &	 116 &	    0.23 &	  8959.02 &	 0/0	 \\ 	 
			 & FLOW         &	 210 &	 210 &	  72 &	    8.16 &	  6966.04 &	 0/0	 \\ 	 
    \hline
    \rule{0pt}{2.6ex}
    \multirow{5}{*}{\begin{minipage}{.9\linewidth}large + region B\end{minipage}}	 & TransQSTP++  &	 309 &	 309 &	 285 &	    4.80 &	  3215.89 &	 206/208	 \\ 	 
			 & TransQSTP+   &	 309 &	 287 &	 262 &	   44.54 &	  3991.13 &	 82/84	 \\ 	 
			 & TransQSTP    &	 309 &	 122 &	 105 &	  281.66 &	 12013.90 &	 0/0	 \\ 	 
			 & QSTP         &	 309 &	  67 &	   3 &	    0.86 &	 21432.33 &	 0/0	 \\ 	 
			 & FLOW         &	 309 &	 290 &	   0 &	   68.50 &	    --    &	 0/0	 \\	 
    \hline\hline
    \end{tabular}
    \label{tab:testset_results_all}
\end{table}

\subsection{Trade-offs between costs and scenic value}\label{sec:tradeoffs}

\begin{figure}
\centering
\begin{subfigure}[b]{.32\linewidth}
\includegraphics[width=\linewidth]{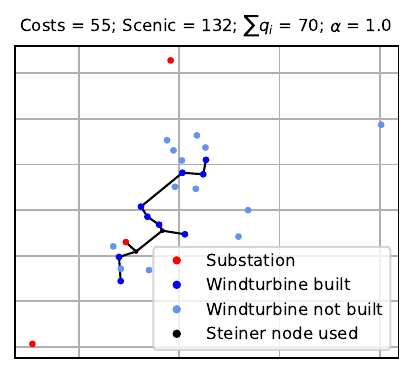}
\caption{}\label{fig:mouse}
\end{subfigure}
\begin{subfigure}[b]{.32\linewidth}
\includegraphics[width=\linewidth]{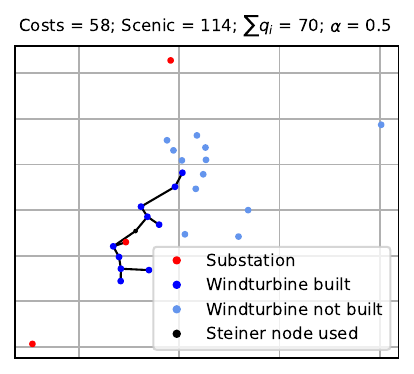}
\caption{}\label{fig:gull}
\end{subfigure}
\begin{subfigure}[b]{.32\linewidth}
\includegraphics[width=\linewidth]{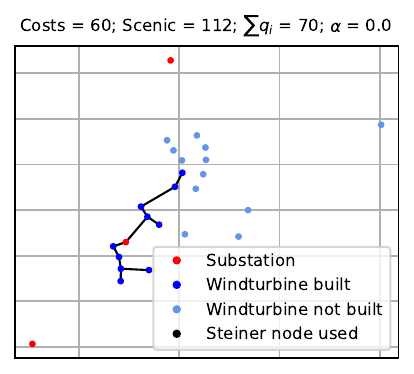}
\caption{}\label{fig:tiger}
\end{subfigure}
\begin{subfigure}[b]{.6\linewidth}
\includegraphics[width=\linewidth]{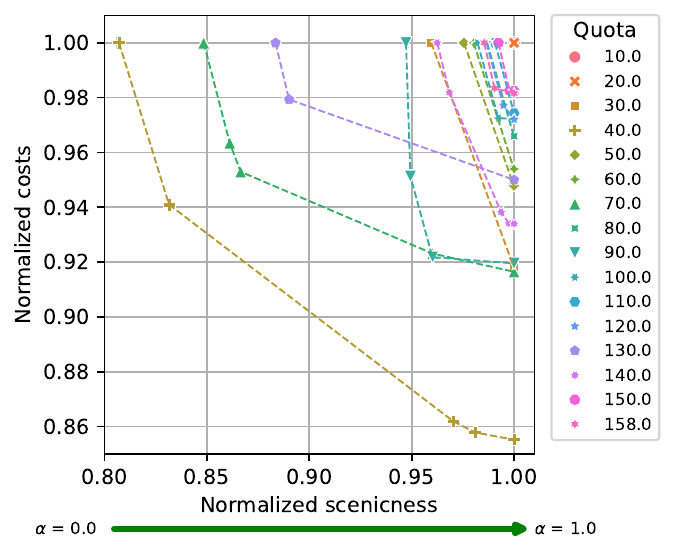}
\caption{}\label{fig:tradeoff_regionA}
\end{subfigure}
\caption{Solution for different $\alpha$ values with a quota of $\quota = 70.0$ of region A (Panels \ref{fig:mouse}-\ref{fig:tiger}). Panel \ref{fig:tradeoff_regionA} shows normalized scenic value versus normalized costs for each quota in region A. The normalization was performed for each curve by dividing the values by the respective maximum cost or scenic value value for the corresponding quota. This means that the curves for different quotas are not comparable to each other, but for each quota the trade-offs between scenic value and cost are shown.}
\label{fig:sol_network_regionA}
\end{figure}

\begin{figure}
\centering
\begin{subfigure}[b]{.49\linewidth}
\includegraphics[width=\linewidth]{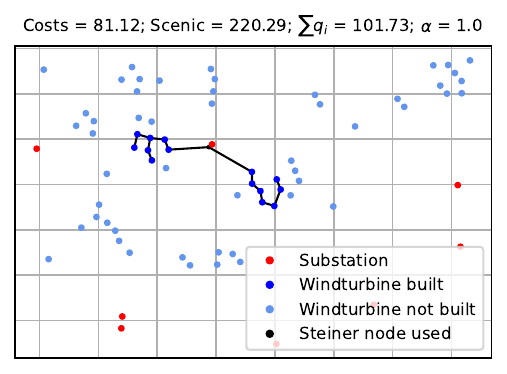}
\caption{}\label{fig:regionB_1.0_100}
\end{subfigure}
\begin{subfigure}[b]{.49\linewidth}
\includegraphics[width=\linewidth]
{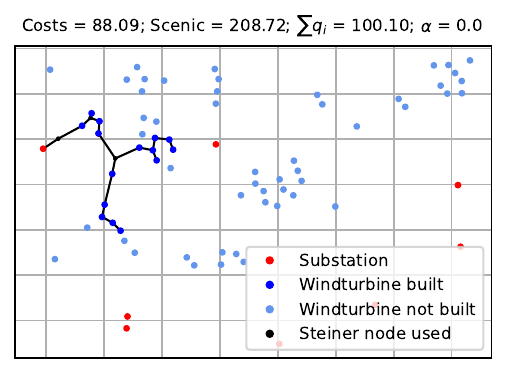}
\caption{}\label{fig:regionB_0.0_100}
\end{subfigure}
\begin{subfigure}[b]{.49\linewidth}
\includegraphics[width=\linewidth]
{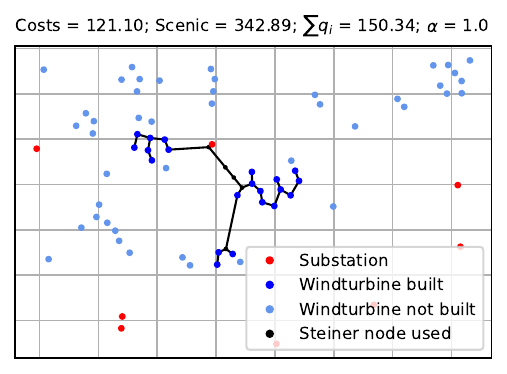}
\caption{}\label{fig:regionB_1.0_150}
\end{subfigure}
\begin{subfigure}[b]{.49\linewidth}
\includegraphics[width=\linewidth]{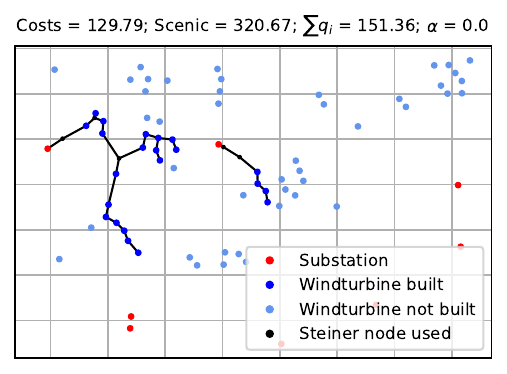}
\caption{}\label{fig:regionB_0.0_150}
\end{subfigure}
\begin{subfigure}[b]{.60\linewidth}
\includegraphics[width=\linewidth]{{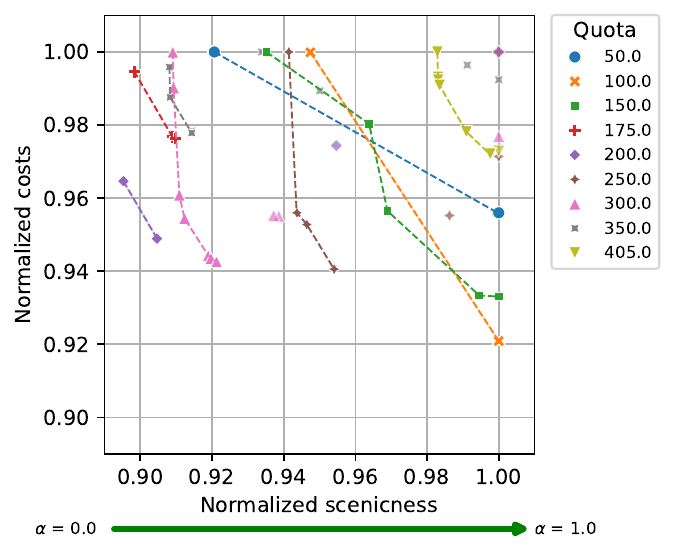}}
\caption{}\label{fig:tradeoff_regionB}
\end{subfigure}
\caption{Solution for a quota of $\quota = 100.0$ with $\alpha=1.0$ (Panel \ref{fig:regionB_1.0_100}) and $\alpha=0.0$ (\ref{fig:regionB_0.0_100}) and for a quota of $\quota = 150.0$ with $\alpha=1.0$ (\ref{fig:regionB_1.0_150}) and $\alpha=0.0$ (\ref{fig:regionB_0.0_150}) of region B. Panel \ref{fig:tradeoff_regionB} shows normalized scenic value versus normalized costs for each quota in region B. The normalization was performed for each curve by dividing the values by the respective maximum cost or scenic value for the corresponding quota. This means that the curves for different quotas are not comparable to each other, but for each quota the trade-offs between scenic value and cost are shown. Points not on their curve are not Pareto optimal, as they are dominated by other solutions. These unconnected points are nonoptimal solutions (optimality gap $>0\%$).}
\label{fig:sol_network_quota_100_regionB}
\end{figure}

The trade-offs between cost and scenic value depend to a large extent on the conditions of the region as well as on the quota, i.e. how much of the maximum possible electricity supply by wind turbines is required (see Figures \ref{fig:sol_network_regionA} and \ref{fig:sol_network_quota_100_regionB}). For low quotas, the solution space offers little flexibility: region A has only one solution for any given $\alpha$ for quotas of 10 and 20 GWh (Figure \ref{fig:tradeoff_regionA}), and region B has only two solutions for the lowest quotas of 50 and 100 GWh (Figure \ref{fig:tradeoff_regionB}). Thus, for the latter problems, no trade-off between cost and scenic value is possible. Either the decision is made to minimize cost while increasing scenic impact by 8\% (quota of 50 GWh) or 5\% (100 GWh), or the focus is on scenic value while increasing costs by 4\% (50 GWh) or 8\% (100 GWh). Likewise, the planning flexibility is low at very high quotas, as in these cases the wind turbines are more or less fixed and the costs and landscape impacts can only be influenced by the connection of the turbines to the substations.

In contrast to the previous insights, planning flexibility is somewhat high for medium quotas. For example, in region A, for a quota of 40 GWh, scenic value can be reduced by almost 20\% or costs by 14\% (Figure \ref{fig:tradeoff_regionA}). A similar trend can be observed for a quota of 70 GWh; for this, Figure \ref{fig:mouse}-\ref{fig:tiger} shows how turbine selection and power networks would vary depending on the weighting of costs and scenic value. The Pareto curves are particularly interesting, however, for the larger region B, such as at quotas of 250 or 300 GWh: in these cases, costs can be reduced by about 5\% with almost unchanged landscape impact (Figure \ref{fig:tradeoff_regionB}).

\subsection{Large-scale siting must include cable routing}\label{sec:Res_Sequentiel_vs_combined}

\begin{figure}[!tbp]
\centering
\begin{subfigure}[t]{.46\linewidth}
\includegraphics[width=\linewidth]{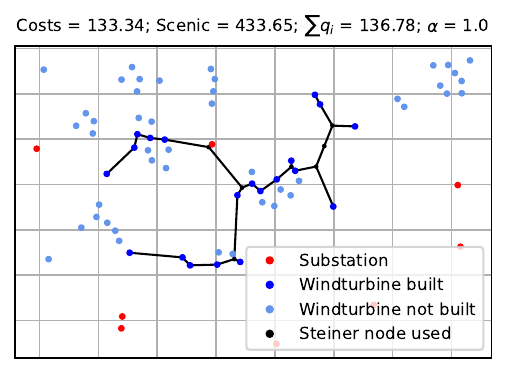}
\caption{Minimize costs: Sequential}\label{fig:D_Sze_200GW_minCosts}
\end{subfigure}
\begin{subfigure}[t]{.46\linewidth}
\includegraphics[width=\linewidth]{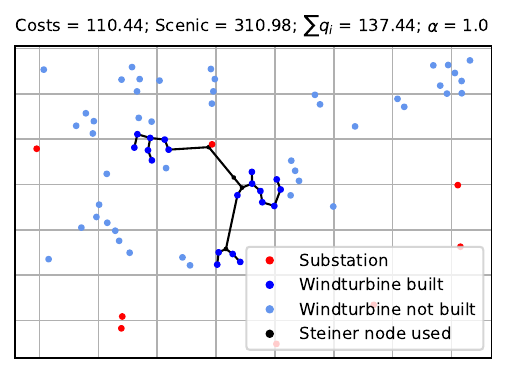}
\caption{Minimize costs: Combined}\label{fig:D_Sze_200GW_minCosts_QSTP}
\end{subfigure}
\begin{subfigure}[t]{.46\linewidth}
\includegraphics[width=\linewidth]{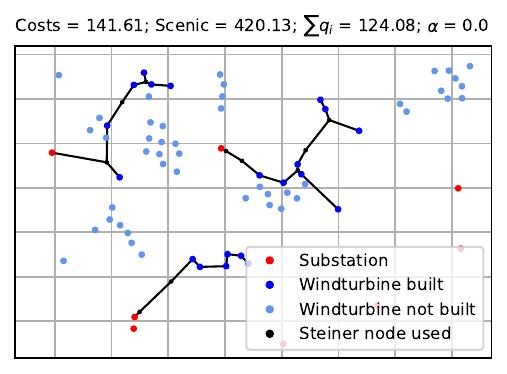}
\caption{Minimize landscape impact: Sequential}\label{fig:D_Sze_200GW_minScenic}
\end{subfigure}
\begin{subfigure}[t]{.46\linewidth}
\includegraphics[width=\linewidth]{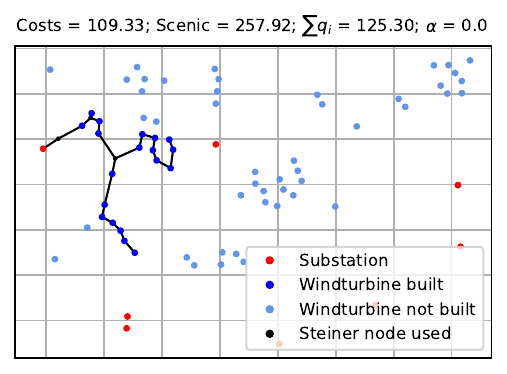}
\caption{Minimize landscape impact: Combined}\label{fig:D_Sze_200GW_minScenic_QSTP}
\end{subfigure}
\caption{Solution for minimizing the costs ($\alpha$ = 1.0) with a quota of around $\quota = 137.0$ in region B (a and b), and for minimizing the scenic value ($\alpha$ = 0.0) with a quota of around $\quota = 124.0$ in region B (c and d). Panels a and c show the optimal turbine locations at the state level for Germany assuming a capacity of 200 GW in 2050 from \textcite{Weinand.29.06.2021}, which did not involve simultaneous planning of the grid connection. Based on the fixed turbine locations, we used our optimization model to connect these turbines to substations. Panels b and d represent the optimal solutions based on our approach with simultaneous turbine and cable routing optimization.}\label{fig:D_Sze_200GW}
\end{figure}

Planning wind turbine expansion on a national level leads to suboptimal solutions in terms of cost and landscape impact if cable routing is neglected. We demonstrate this by comparing our approach with the scenarios from \textcite{Weinand.29.06.2021}, who identified the optimal turbine locations at the national level for Germany assuming a capacity of 200 GW in the year 2050 (Figure \ref{fig:D_Sze_200GW}). In \textcite{Weinand.29.06.2021}, only the costs and scenic value of the turbine sites were taken into account, but not for the cable routing, as a simultaneous optimization of both factors was hardly computationally practicable for such a large region. In order to comparatively assess the turbine locations in \textcite{Weinand.29.06.2021} using our approach, we first connected the turbines, found in \textcite{Weinand.29.06.2021} for region B, to the substations while minimizing costs or scenic value. This is consistent with the methodology developed in this study, but with turbine sites chosen in advance. Afterwards, for both scenarios, i.e. only costs and only scenic value, we solve the QSTP for the quota values given by \textcite{Weinand.29.06.2021} for this region, i.e., $\quota=136.78$ and $\quota=124.08$ for only costs and only scenic value, respectively.
When minimizing costs, selecting turbine sites without simultaneously considering the cable costs can result in 21\% higher total costs, i.e., costs for turbines and cable routing, than simultaneously optimizing these costs (Figure \ref{fig:D_Sze_200GW_minCosts} and Figure \ref{fig:D_Sze_200GW_minCosts_QSTP}). The same is true when minimizing only the landscape impact. The landscape impact is 38.6\% lower with the simultaneous optimization of turbines and the cable routing compared to the sequential approach (Figure \ref{fig:D_Sze_200GW_minScenic} \& Figure \ref{fig:D_Sze_200GW_minScenic_QSTP}). Interestingly, the shorter cable routing chosen in our novel approach would always also lead to lower values for the criterion that was not been chosen as the objective. For example, if scenic value were to be minimized for a quota of Q = 137 GWh (as for the cost minimzations in Figure \ref{fig:D_Sze_200GW_minCosts} \& Figure \ref{fig:D_Sze_200GW_minCosts_QSTP}), the overall costs would still be 10\% lower than for the approach with fixed turbine siting and subsequent cost-optimized cable routing. 

In the future, it is vital that recent approaches in the articles by 
\textcite{Weinand.29.06.2021}, \textcite{Lehmann.2021b}, \textcite{Tafarte.2021}, and \textcite{SPIELHOFER2023220} on the optimal siting of turbines on a national level attempt to simultaneously include the cable routing. For this purpose, approaches should be developed to make the methodology used in this study applicable to very large regions, such as entire countries.

\section{Discussion and conclusions} \label{sec:Discussion}

Current methods for planning onshore wind farms and connecting them to the power network do not always meet minimum cost solutions or account for social and environmental factors. This paper proposes a new approach for optimizing turbine locations and cable routing using a Quota Steiner tree problem (QSTP). We present a novel transformation of the already known directed cut formulation of the QSTP. Even though the LP-relaxations of the transformed and original formulation are proven to be equivalent, the transformed formulation proves to be more effective in practice. By using shortest path reduction and the shortest path heuristic, we have advanced the state-of-the-art STP-solver \scipjack in the context of QSTP, which outperforms standard MIP solvers and is capable of solving large-scale problems with up to at least 1.3 million edges. Our case studies, conducted in selected regions of Germany, demonstrate significant trade-offs between cost and landscape impacts when planning onshore wind farms. However, we also show that small reductions in one objective criterion can yield significant improvements in the other. Additionally, we illustrate that simultaneously optimizing grid integration with turbine siting is essential for avoiding excessive costs or landscape impacts in the course of wind farm projects. 

In this study, our focus was on presenting the applicability and utility of the methodology. For real case studies, further important criteria beyond costs and landscape impacts should be considered in the future. For example, disamenities for the local population should be incorporated \parencite{Lehmann.2021b, Tafarte.2021, Weinand.historic}, as well as the environmental impacts of turbines, e.g., through bird strikes \parencite{McKenna.2022} and land use competition with other renewable energy sources \parencite{MCKENNA2022123754}. Furthermore, we used different weightings in our multi-objective planning without knowing the exact preferred weightings between the two objectives costs and scenic value. In addition, equating the landscape impact of a wind turbine with that of a 1 km electricity network is a strong simplification. However, previous research has demonstrated that finding universally valid weights between target criteria in wind farm planning is nearly impossible, even for experts \parencite{Lehmann.2021}. In regional planning, such as that discussed in this paper, local stakeholders could be consulted in multi-criteria decision approaches to determine appropriate weightings in the future \parencite{MCKENNA20181092}. By applying other multi-objective methods, such as the $\varepsilon$-constraint methods, one could additionally diversify the set of non-dominated solutions.

Furthermore, we select the Steiner points independently of the geographical conditions. Our problem could be improved by using a digital elevation model to determine true distances and by taking obstacles into account for which the cables cannot be deployed (see also \textcite{Fischetti.2018}). If only Euclidean distances between potential turbines are considered, instead of defining a discrete set of Steiner points, one could also approach the problem as an extension of the Euclidean Steiner tree problem with obstacles; for an overview see, e.g., \textcite{Brazilbook}. However, since we do not consider Euclidean distance, state-of-the-art algorithms, e.g., implemented in the \texttt{GeoSteiner} software by \textcite{Juhl2018}, do not work in our case. In addition, for the set of the larger regions, we were not able to find the optimal solution for all instances. Especially when transferring our methodology to federal states or countries, further methodological advancements will be necessary to make this computable. For example, reduction techniques are one of the most essential and effective features in solving STP-related problems \parencite{ljubic2021} and \scipjack provides many strong reduction techniques, see e.g., \textcite{RehfeldtKoch2023, rehfeldt2021a, rehfeldt2019}. Therefore, identifying efficient reduction techniques could significantly improve the solution process, however, these have yet to be studied in the context of QSTPs. Additionally, we used a shortest path heuristic to find primal solutions in the branch-and-cut algorithm. In the future, developing new primal and dual heuristics in the context of QSTP could be a promising way to further improve our approach. For instance, exploration of the idea of \textcite{leitner2018} who introduced a dual-ascent algorithm in the context of PCSTPs could be worthwhile. When investigating larger regions and wind farms, future studies should also incorporate the meshing of power grids for increased security of supply as well as the (remaining) capacities of substations. The existing substations we have considered may not always have sufficient remaining capacity to ensure the grid connection of further wind farms and new substations would have to be installed. Also, cables are subject to finite power capacities, which are not included in the model, but are inevitable for the realization of particular turbine sites, see, e.g., \textcite{Fischetti.2018}. Future studies should introduce additional constraints representing other relevant technical details as well, such that our approach becomes more applicable for single wind farm planning. For example, including a minimum distance between the turbines; Or adapting the quota constraint to model the wake effect explicitly in the optimization. Then our approach could be combined by the cuts introduced in \textcite{fischetti2022}. Although our approach is also applicable for planning other renewable energy plants like solar photovoltaic installations, applications of the algorithm beyond energy system analysis are also conceivable: for example, cable routing optimization approaches for offshore wind farms \parencite{Fischetti.2018, Fischetti.2021a} were recently applied to determine safe distancing during the COVID-19 pandemic \parencite{Fischetti.2021b}.  
 
{\singlespacing
\section*{Acknowledgements}
The work for this article has been conducted in the Research Campus MODAL funded by the German Federal Ministry of Education and Research (BMBF) (fund numbers 05M14ZAM, 05M20ZBM). Furthermore, this work was supported by the Helmholtz Association under the program "Energy System Design".
\singlespacing 
\section*{CRediT statement}
Conceptualization: J.W., J.P., D.R.; data curation: J.W., J.P., C.S.; formal analysis: J.P., J.W.; investigation: J.P., J.W.; methodology: J.P., D.R., J.W.; software: J.P., D.R.; validation: J.P., J.W.; visualization: J.P., J.W.; writing – original draft: J.P., J.W.; writing – review and editing: J.P., J.W., D.R. 
}

{\singlespacing
\printbibliography
}

\end{document}